\documentclass[12pt, twoside]{amsart}
\usepackage{amsmath,amssymb,amscd,amsxtra,latexsym,
graphicx,epsfig,euscript,amsthm,mathscinet}
\usepackage{graphicx}
\usepackage{float}
\usepackage{a4wide}
\usepackage{amsmath,amsfonts,amsthm,amsopn,amssymb,enumitem}
\usepackage{palatino}
\usepackage[backend=bibtex,style=numeric,sorting=nyt,
giveninits=true,isbn=false,url=false,doi=false]{biblatex}
\renewbibmacro{in:}{\ifentrytype{article}{}{\printtext{\bibstring{in}\intitlepunct}}}
\DeclareFieldFormat*{title}{#1}
\usepackage[colorlinks]{hyperref}
\AtBeginDocument{%
  \hypersetup{%
    linkcolor=blue,%
    citecolor=green,%
  }%
}
\usepackage{cleveref}
\crefname{equation}{}{}
\crefname{figure}{}{}

\theoremstyle{plain}
\newtheorem{theorem}{Theorem}[section]
\newtheorem*{theorem*}{Theorem}
\newtheorem{lemma}[theorem]{Lemma}
\newtheorem{proposition}[theorem]{Proposition}
\newtheorem{corollary}[theorem]{Corollary}

\theoremstyle{definition}
\newtheorem*{remarks*}{Remark}
\newtheorem{remark}[theorem]{Remark}

\newtheorem*{example*}{Example}
\newtheorem*{examples*}{Examples}

\newtheorem*{definition*}{Definition}

\newcommand{\proofend}{\hspace*{\fill} $\Box$\\}

\def\1{\:\!}
\def\2{\;\!}

\def\Diffc0{\operatorname{Diff^c_0}}

\def\Sympc0{\operatorname{Symp^c_0}}

\def\pp{\partial}

\def\.{\mskip1mu}
\def\?{\mskip-1mu}

\def\1{\:\!}
\def\2{\;\!}

\def\Diffc0{\operatorname{Diff^c_0}}

\def\Sympc0{\operatorname{Symp^c_0}}

\def\pp{\partial}

\def\.{\mskip1mu}
\def\?{\mskip-1mu}

\def\proof{\noindent {\it Proof. \;}}

\def\R{\mathbb{R}}

\def\tf{\tilde{f}}
\def\tF{\tilde{F}}

\begin{document}

\title[Rigidity results for the capillary overdetermined problem]{Rigidity results for the capillary overdetermined problem}

\author{Yuanyuan Lian}
\address{(Y.~Lian)
Departamento de An\'alisis matem\'atico, Universidad de Granada,
Campus Fuen-tenueva, 18071 Granada, Spain} \email{lianyuanyuan.hthk@gmail.com; yuanyuanlian@correo.ugr.es}

\author{Pieralberto Sicbaldi}
\address{(P.~Sicbaldi)
Departamento de An\'alisis matem\'atico and IMAG, Universidad de Granada,
Campus Fuentenueva, 18071 Granada, Spain \&
Aix Marseille Universit\'e - CNRS, Centrale Marseille - I2M, UMR 7373, 13453 Marseille, France}
\email{pieralberto@ugr.es}


\maketitle
\begin{abstract}
In this paper we obtain rigidity results for bounded positive solutions of the general capillary overdetermined problem
\begin{equation}\label{cap_intro}
  \left\{\begin{array} {ll}
\mathrm{div} \left(\frac{\nabla u}{\sqrt{1+|\nabla u|^2}}\right) + f(u) = 0 & \mbox{in }\; \Omega,\\
               u= 0 & \mbox{on }\; \pp \Omega, \\
\partial_{\nu} u=\kappa &\mbox{on }\; \pp \Omega,
\end{array}\right.
\end{equation}
where $f$ is a given $C^1$ function in $\mathbb{R}$, $\nu$ is the exterior unit normal, $\kappa$ is a constant and $\Omega \subset \mathbb{R}^n$ is a $C^1$ domain. Our main theorem states that if $n=2, \kappa\neq 0$, $\pp \Omega$ is unbounded and connected, $|\nabla u|$ is bounded and there exists a nonpositive primitive $F$ of $f$ such that $F(0)\geq \left(1+\kappa^2\right)^{-\frac12} -1$, then $\Omega$ must be a half-plane and $u$ is a parallel solution. In other words, under our assumptions, if a capillary graph has the property that its mean curvature depends only on the height, then it is the graph of a one dimensional function. We also prove the boundedness of the gradient of solutions of \eqref{cap_intro} when $f'(u) <0$. Moreover we study a Modica type estimate for the overdetermined problem \eqref{cap_intro} that allows us to prove that, unless $\Omega$ is a half-space, the mean curvature of $\pp \Omega$ is strictly negative under the assumption that $\kappa\neq 0$ and there exists a nonpositive primitive $F$ of $f$ such that $F(0)\geq \left(1+\kappa^2\right)^{-\frac12} -1$. Our results have an interesting physical application to the classical capillary overdetermined problem, i.e., the case where $f$ is linear.
\end{abstract}

\section{Introduction and statement of the main results}\label{sec.Intro}

The classical capillary phenomenon is given by a liquid rising in a straight cylindrical tube dipped into a large (mathematically infinite) reservoir. Let $\Omega \subset \R^2$ be the section of the tube and $u(x,y)$ be the height of the upper surface of the liquid (relative to the reservoir). Then $u$ is proportional to the mean curvature of the upper surface (see \cite[Chapter 1.2]{MR816345}). In addition, the surface of the liquid intersects the wall of the tube with a constant angle $\theta$, called {\it wetting angle} (see \cite[Chapter 1.3]{MR816345}). The wetting angle can be less than $\pi/2$ (for example in the case of water) or bigger than $\pi/2$ (for example for mercury). As usual, we assume the former case in this paper. The later case can be transformed to the former case by considering $v=-u$ (see \cite[Chapter 1.9 on Page 13]{MR816345} and \cite[Section 3 on Page 477]{MR607989}). Therefore, $u$ satisfies (see \cite{MR670443, MR670442,MR607989}):
\begin{equation}\label{equil}
  \left\{\begin{array} {ll}
\mathrm{div} \left(\frac{\nabla u}{\sqrt{1+|\nabla u|^2}}\right) - b\, u = 0 & \mbox{in }\; \Omega,\\
                \partial_{\nu} u = \cos \theta\, \sqrt{1+ |\nabla u|^2} & \mbox{on }\; \pp \Omega,
\end{array}\right.
\end{equation}
where $\nu$ is the exterior unit normal on $\pp \Omega$. Here, $b = \frac{(\rho-\rho_0)\, g}{\sigma}>0$, where $\rho, \rho_0, \sigma, g$ are respectively the density of the liquid, the density of the gas outside the liquid, the surface tension of the liquid and the gravity.

Now, if we assume that the liquid raises the same height on the wall of the tube, then $u$ satisfies
\begin{equation}\label{capil-0}
  \begin{cases}
  \mathrm{div} \left(\frac{\nabla u}{\sqrt{1+|\nabla u|^2}}\right) - b u =0& \mbox{in $\Omega$, }\\
  u=c_h &\mbox{on $\partial\Omega$, }\\
  \partial_{\nu} u=\kappa_0 &\mbox{on $\partial\Omega$\,, }
  \end{cases}
\end{equation}
where $c_h>0$ and $\kappa_0=\cot \theta>0$ are constants. In physics it was very well known since many years that, unless the wetting angle is equal to $\pi/2$ the only case when the liquid raises the same height on the wall of the tube is when $\Omega$ is a disk (Jurin's law, 1718), but mathematically this rigidity result is due to Serrin \cite{MR333220}. Let us explain a bit more. Since $\Omega$ is bounded, by the comparison principle, $0<u< c_h$ in $\Omega$. Let $v=c_h-u$ and then $v$ is a solution of
\begin{equation}\label{capil}
  \begin{cases}
  \mathrm{div} \left(\frac{\nabla v}{\sqrt{1+|\nabla v|^2}}\right) - b v +bc_h=0& \mbox{in $\Omega$, }\\
  v>0 &\mbox{in $\Omega$, }\\
  v=0 &\mbox{on $\partial\Omega$, }\\
  \partial_{\nu} v=\kappa &\mbox{on $\partial\Omega$\,, }
  \end{cases}
\end{equation}
where $\kappa=-\kappa_0<0$. We refer to \eqref{capil} as the {\it classical capillary overdetermined problem}. Serrin \cite[Theorem 2]{MR333220} proved that if $\Omega$ is a bounded regular domain of $\R^n$ ($C^1$ is enough) and $f$ is a $C^1$ function such that there exists a solution to \cref{capil}, then $\Omega$ must be a ball and $u$ must be radial. In fact, Serrin proved the symmetry for a more general overdetermined elliptic problem
\begin{equation}\label{serrin}
  \left\{\begin{array} {ll}
Q(u, |\nabla u|) + f(u, |\nabla u|) = 0 & \mbox{in }\; \Omega,\\
u> 0 & \mbox{in }\; \Omega, \\
               u= 0 & \mbox{on }\; \pp \Omega, \\
\partial_{\nu} u=\kappa &\mbox{on }\; \pp \Omega,
\end{array}\right.
\end{equation}
where $Q(u)$ is the operator given by
\[
Q(u, |\nabla u|): = a(u,|\nabla u|) \Delta u + h(u,|\nabla u|)\, u_i u_j u_{ij}
\]
for some two variables $C^1$ functions $a, h, f$.

\medskip

Another interesting case of the classical capillary problem is when we consider a large reservoir of a liquid and we dip a finite number, say $m$, of cylindrical solids with section $G_i$ ($i=1, ..., m$). If $u$ represents the surface of the liquid and we assume that the liquid raises the same height on the walls of the cylinders, then we will get again problem \eqref{capil} where now $\Omega = \R^n \backslash \cup G_i$. By some results on overdetermined elliptic problem in the form \eqref{serrin} in exterior domains due to Reichel \cite{MR1463801} and Sirakov \cite{ MR1808026}), the only possibility is that $m=1$ and $G_1$ is a disk.

\medskip

\begin{figure}[!ht]
  \centering
  \includegraphics[width=10cm]{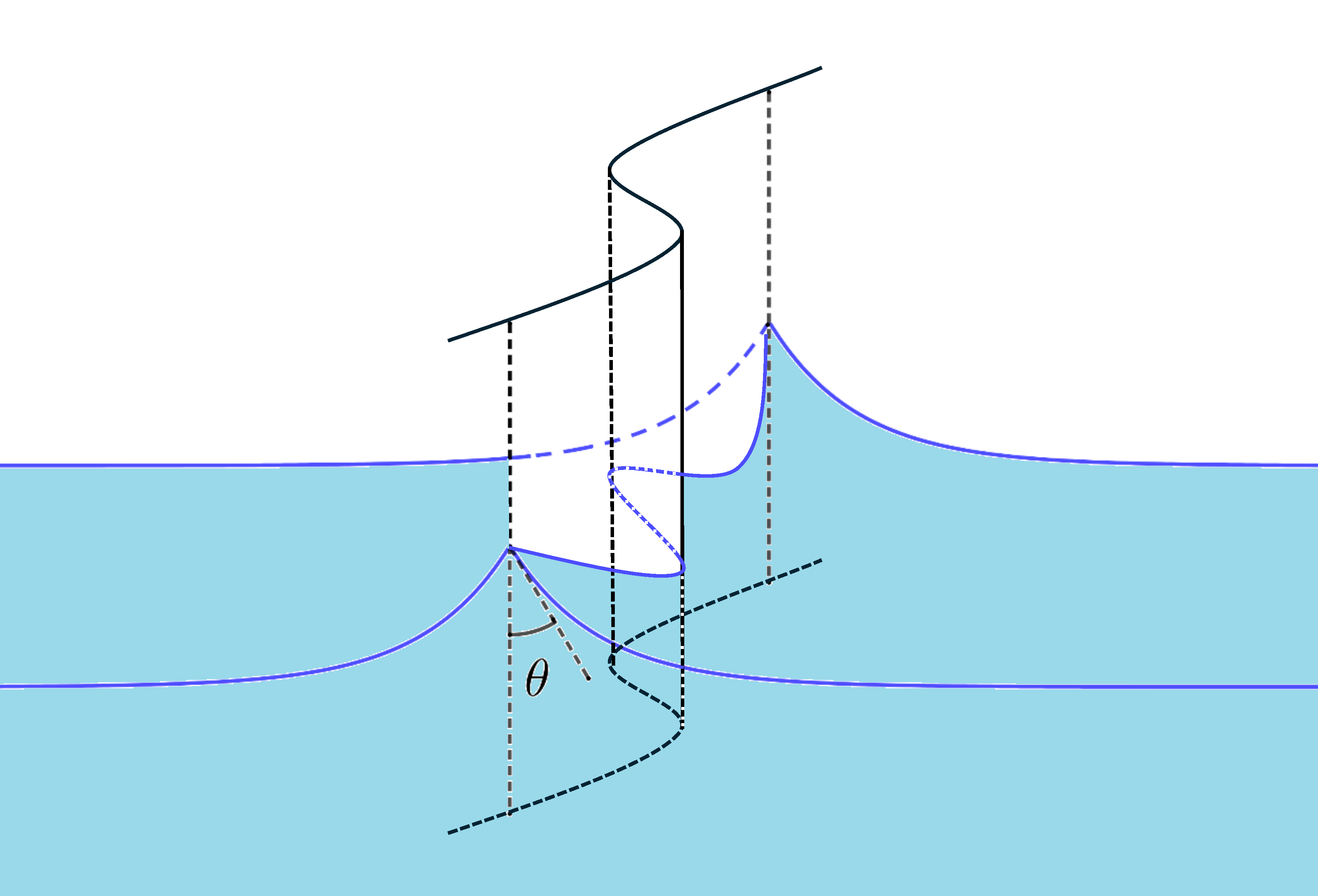}
  \caption{The classical capillary phenomenon inserting a vertical plate in a reservoir of liquid.}\label{fig1}
\end{figure}

In this paper we will consider the third interesting case. As before, we consider a liquid in a large reservoir  and we dip in it a large (mathematically infinite) cylindrical plate, i.e., a rigid cylinder whose section is a regular open unbounded curve. Clearly the surface of the liquid will assume a certain nonconstant shape intersecting the plate with a constant wetting angle $\theta$ on both sides of the plate, see FIGURE \Cref{fig1}. Without loss of generality, here we just consider the capillary phenomenon only on one of the two connected components of the space separated by the plate. The function $u(x,y)$ that represents the surface of the liquid is then defined in a domain $\Omega$ diffeomorphic to a half-plane and satisfies \eqref{equil}. Clearly, the height of the liquid on the plate in general is not constant, so we can pose the following:

\medskip

\noindent\textbf{Question 1:} For which possible shapes of the plate the liquid will raise the same height on the plate?


\begin{figure}[!ht]
  \centering
  \includegraphics[width=10cm]{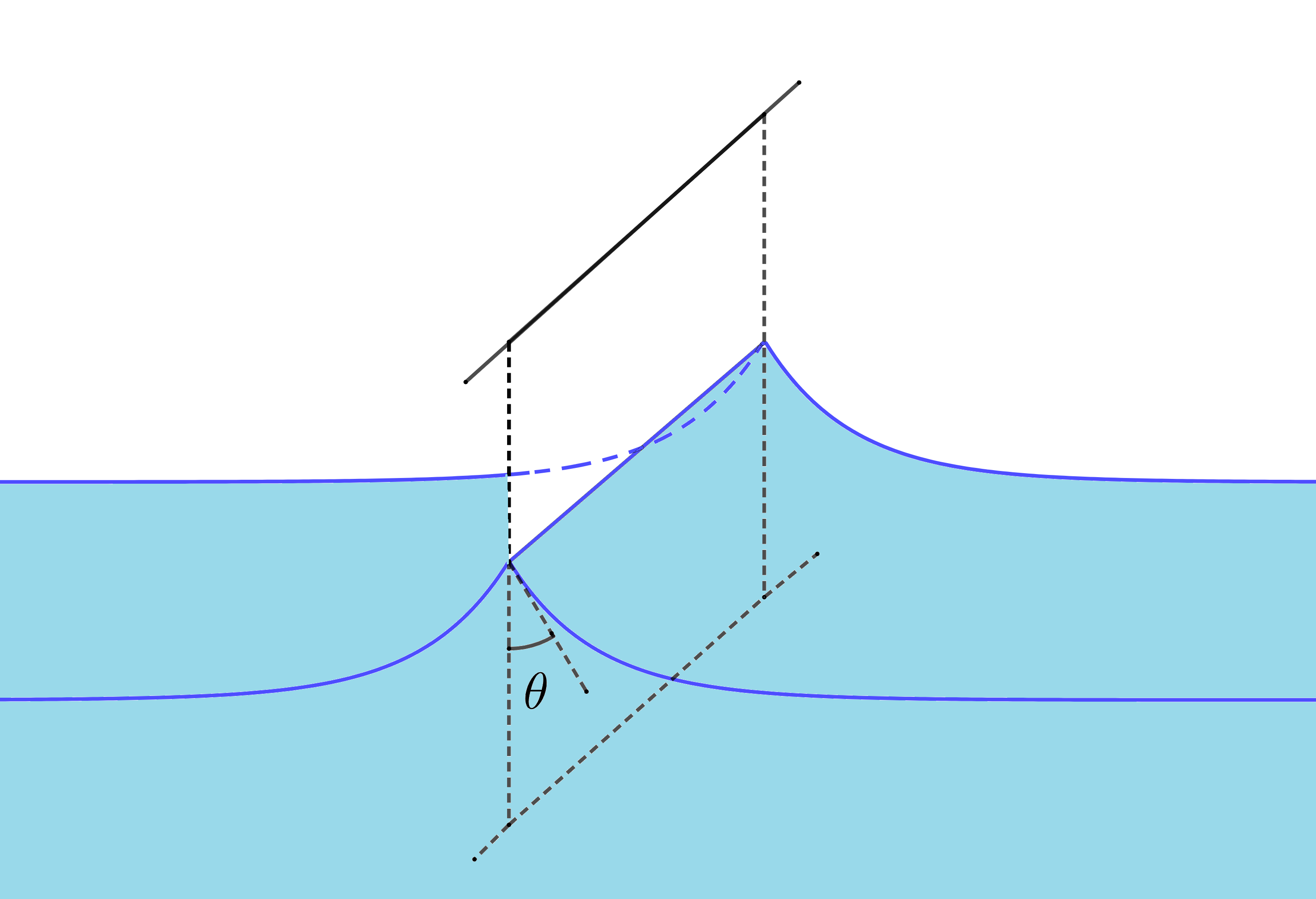}
  \caption{The case where the plate is planar.}\label{fig2}
\end{figure}

Up to our knowledge this is an open question, and in this paper we will give an answer to it. In physics it is well known that the liquid rises the same height when the plate is planar (see FIGURE \Cref{fig2}). This fact is tightly linked to the Wilhelmy method to measure the surface tension of a liquid (see \cite{Wilhelmy1863, MILLER2001383, TIAB2016319, RAPP2017453}). From the mathematical point of view, the question is: is the half-plane the only domain $\Omega \subset \R^2$ with boundary unbounded and connected where a solution of \cref{capil-0} can exist? In this paper we prove that the answer to this question is affirmative (see \Cref{cor_cap}).

It is more convenient to consider a positive solution in rigidity problems. By the comparison principle for unbounded domains (see \cite{MR991023}), $0<u< c_h$ in $\Omega$. As before, let $v=c_h-u$ and then $v$ is a solution of \cref{capil}. In fact, we will consider a more general overdetermined problem, i.e.,
\begin{equation}\label{e.R}
  \begin{cases}
  \mathrm{div} \left(\frac{\nabla u}{\sqrt{1+|\nabla u|^2}}\right)+f(u) =0& \mbox{in $\Omega$, }\\
  u>0 &\mbox{in $\Omega$, }\\
  u=0 &\mbox{on $\partial\Omega$, }\\
  \partial_{\nu} u=\kappa &\mbox{on $\partial\Omega$, }
  \end{cases}
\end{equation}
where $f$ is a given $C^1$ function in $\mathbb{R}$. Such problem arises again in the capillary phenomenon, when the liquid is subject to an external force different from the gravity. See for example \cite{MR2335074} where cases of nonlinear functions $f(u)$ are physically considered. We refer to \eqref{e.R} as the {\it general capillary overdetermined problem.} Observe that in the case $f\equiv H$, where $H$ is a constant, we have the usual geometric case of capillary constant mean curvature graphs. For other functions $f$, The \eqref{e.R} means that we have the graph of a function $u$ defined in a domain $\Omega$, that intersects the horizontal plane with constant angle on $\partial \Omega$ and such that its mean curvature depends, in any point, only on the height of the graph.

\medskip

Overdetermined elliptic problems in unbounded domains have been mainly considered in the form \eqref{serrin} with $Q(u, |\nabla u|)=\Delta u$ and $f$ just depending on $u$, i.e.,
\begin{equation}\label{bcnn}
  \left\{\begin{array} {ll}
\Delta u + f(u) = 0 & \mbox{in }\; \Omega,\\
u> 0 & \mbox{in }\; \Omega, \\
               u= 0 & \mbox{on }\; \pp \Omega, \\
\partial_{\nu} u=\kappa &\mbox{on }\; \pp \Omega\,.
\end{array}\right.
\end{equation}
In \cite{MR1470317} Berestycki, Caffarelli and Nirenberg obtained rigidity results in the class of domain that are epigraphs, and proposed the very well known following conjecture, that has represented the starting point for all the next research on overdetermined elliptic problems in unbounded domains:

\medskip

{\bf Conjecture BCN}: If $\mathbb{R}^n \backslash
\overline{\Omega}$ is connected, then the existence of a bounded
solution to problem \eqref{bcnn} implies that $\Omega$ is either
a ball, a half-space, a generalized cylinder $B^k \times
\mathbb{R}^{n-k}$ ($B^k$ is a ball in $\mathbb{R}^k$), or the
complement of one of them.

\medskip

It turns out that the conjecture is not true. The first counterexample has been given by the second author in \cite{MR2592974} (see also \cite{MR2854185}), where problem \eqref{bcnn} for the PDE $\Delta u + \lambda\, u$ has been solved in domains $\Omega$ in $\R^n$, for $n \geq 3$, given by periodic rotationally symmetric perturbations of a straight full cylinder. Other counterexamples have been constructed in \cite{MR3417183, MR3590664, MR4046014, MR4484836}. See also the survey \cite{MR4559541}. These counterexamples showed that the geometry of solutions to overdetermined elliptic problems is richer that originally expected. Nevertheless in dimension $n=2$ Ros, Ruiz and the second author proved in the following rigidity result (\cite{MR3666566}):

\medskip

{\bf Theorem RRS}: Suppose that $\Omega$ is a domain of $\R^2$ with boundary unbounded and connected, and assume that $u$ is a bounded solution to \eqref{bcnn} for some Lipschitz function $f$, with $\kappa \neq 0$. Then $\Omega$ must be a half-plane and $u$ is parallel.

\medskip

In general, we say that a solution $u$ is parallel (or one-dimensional) if there exist $x_0,a \in \R^n$ and a function $g: [0, +\infty) \rightarrow\mathbb{R}$ such that
\[
u(x)=g(a\cdot (x-x_0)), \ \ x \in  \{x \in \R^n:  \ a \cdot (x-x_0) >0 \}.
\]

%

\medskip

Overdetermined elliptic problems of the form \eqref{e.R} in unbounded domains, up to our knowledge, have not been studied in a systematic way and few results are available. A priori, the question raised in the original BCN conjecture could be addressed to such overdetermined elliptic problems, and would be interesting to prove similar rigidity results and to construct new nontrivial solutions. This question is in fact suggested in \cite{MR4271788}, where the authors consider problem \eqref{e.R} when $f \equiv H$, where $H$ is a constant. In this paper we aim to generalize the Theorem RRS to the general capillary overdetermined problem.
Our main result is the following:
\begin{theorem} \label{Tmain}
Let $\Omega \subset \mathbb{R}^2$ be a $C^{1}$ domain whose boundary is unbounded and connected. Let $f\in C^1$ and $u\in C^3$ be a solution of \cref{e.R} with
\begin{equation}\label{e1.1}
u\leq C,\quad |\nabla u|\leq C
\end{equation}
for some constant $C$. Suppose that $\kappa \neq 0$ and there exists a nonpositive primitive $F\in C^2(\mathbb{R})$ of $f$ such that
\begin{equation}\label{e.F}
F(0) \geq \left(1+\kappa^{2}\right)^{-\frac{1}{2}} -1.
\end{equation}
Then $\Omega$ is a half-plane and $u$ is parallel.
\end{theorem}

If there exists a bounded parallel solution to problem \eqref{e.R} in the half-plane $\mathbb{R}_{+}^2$, then there exists a nonpositive primitive $F$ of $f$ such that (see \Cref{th1.4})
\begin{equation*}
  F(0) = \left(1+\kappa^{2}\right)^{-\frac{1}{2}} -1.
\end{equation*}
Hence, an alternative of \Cref{Tmain} is
\begin{corollary}\label{co1.1}
Let $\Omega \subset \mathbb{R}^2$ be a $C^{1}$ domain whose boundary is unbounded and connected. Let $f\in C^1$ and $u\in C^3$ be a solution of \cref{e.R} with
\begin{equation}\label{e1.1}
u\leq C,\quad |\nabla u|\leq C
\end{equation}
for some constant $C$. Suppose that $\kappa\neq 0$ and there exists a bounded parallel solution to problem \eqref{e.R} in the half-plane $\mathbb{R}_{+}^2$. Then $\Omega$ is a half-plane and $u$ is parallel.
\end{corollary}

\begin{remark}\label{re1.2}
It is tempting to conjecture that if $f$ is such that condition \cref{e.F} is not satisfied for any nonpositive primitive $F$ of $f$, then problem \cref{e.R} does not have any solution in domains of $\R^2$ with boundary unbounded and connected. This is an open problem.
\end{remark}

\begin{remark}\label{re1.1}
We ask the solution $u$ is bounded and also has bounded gradient. This is a difference with respect to the original BCN conjecture. For the PDE with the Laplace operator (problem \cref{bcnn}), the interior gradient bound can be derived directly from the classical Schauder regularity (see \cite[Theorem 4.5]{MR1814364}). The boundary gradient bound is more complicated, and can be obtained using curvature estimates of $\partial \Omega$ (see \cite[Chapter 3]{MR3666566}) or the method in \cite[Lemma 2.2]{Ruiz_Sicbaldi_Wu24}. For problem \cref{e.R}, it is well-known that one cannot get $|\nabla u|\leq C$ for a general $f$. Even for the interior gradient bound, one usually requires $f'\leq 0$ (see \cite{MR843597, MR1617971}). Therefore, we assume $|\nabla u|\leq C$ in \Cref{Tmain}. In addition, we notice that the PDE in \cref{e.R} is uniformly elliptic if $\nabla u$ is bounded, and uniformly elliptic PDEs have a rich regularity theory, which is very important for the rigidity.
\end{remark}

\medskip

Nevertheless we are able to prove the following regularity result, that holds in any dimension $n$:
\begin{theorem}\label{prop2.1}
Let $\Omega \subset \mathbb{R}^n$ be a $C^{1}$ domain, $u\in C^3$ be a bounded solution of \cref{e.R} for some function $f \in C^1$ and some constant $\kappa$. 
Assume that $f'(u)\leq 0$. Then
\begin{equation}\label{e2.13}
  ||\nabla u||_{L^{\infty}(\Omega)}\leq C,
\end{equation}
where $C$ depends only on $n, \kappa, \|u\|_{L^{\infty}(\Omega)}, \|f\|_{L^{\infty}}$.
\end{theorem}

\medskip

By applying \Cref{Tmain} and \Cref{prop2.1} together, we obtain a complete answer to \textbf{Question 1} about the classical capillary problem we stated before, i.e.,
\begin{theorem}\label{cor_cap}
Let $\Omega \subset \mathbb{R}^2$ be a domain diffeomorphic to a half-plane and $u$ be a solution of \cref{capil-0} with $\kappa_0 \neq 0$ (i.e., the wetting angle is not $\pi/2$). 
Then $\Omega$ is a half-plane, i.e., the infinite plate must be planar (see FIGURE \Cref{fig2}).
\end{theorem}


\medskip


The proof of \Cref{Tmain} is partially inspired by the strategy of the proof of Theorem RRS. Nevertheless, a large part of such work cannot be adapted to our case because the two overdetermined problems are different. In fact, in \cite{MR3666566} the authors use a technique inspired by the proof of the De Giorgi conjecture in dimension 3 (\cite{MR1775735}) to derive that if the solution $u$ of the overdetermined problem is monotone in one variable then the domain must be a half-plane (see \cite[Section 4]{MR3666566}). So an important part of the proof of Theorem RRS is devoted to prove that the solution $u$ is monotone in one variable. Unfortunately, this strategy depends heavily on the special structure of the Laplace operator and seems not applicable to \cref{e.R}. So, in this paper we adopt a completely different strategy, based on the $P$-function. In \Cref{sec2} we develop a Modica type estimate for problem \eqref{e.R}, inspired by the classical Modica estimate \cite{MR803255} and we get from it a rigidity result for overdetermined problem in the form \eqref{e.R} valid in any dimension $n$. More precisely, we will prove the following:
\begin{theorem}\label{coro}
Let $\Omega \subset\mathbb{R}^{n}$ be a $C^1$ domain, $f$ a $C^1$ function and $\kappa \neq 0$. Let $u$ be a $C^3$ solution of \eqref{e.R} such that $u\leq C$ and $|\nabla u|\leq C$ for some constant $C$. If there exists a nonpositive primitive $F$ of $f$ such that
\begin{equation} \label{cond1}
 F(0) \geq \left(1+\kappa^{2}\right)^{-\frac{1}{2}}  -1 ,
\end{equation}
then either $\Omega$ is a half-space and $u$ is parallel, or $H(p) < 0$ for any $p \in \partial\Omega$, where $H(p)$ is the mean curvature of $\partial \Omega$ at the point $p$.
\end{theorem}

\medskip

We will obtain the previous \Cref{coro} as a corollary of a more general result, that represents the Modica type estimate for the elliptic overdetermined problem \eqref{e.R}. Such estimate is the counterpart of the estimate obtained by Ruiz, Wu and the second author in \cite{Ruiz_Sicbaldi_Wu24} for the overdetermined problem \eqref{bcnn}. We believe that such estimate itself is an interesting result for the study of overdetermined problems in \eqref{e.R}, and we remind the reader to \Cref{sec2} for its precise statement (\Cref{Th.P1}).

\medskip

From \Cref{coro}, an immediately corollary can be obtained.
\begin{corollary}\label{coro2}
Suppose that $u$ is a $C^3$ solution of \eqref{e.R} in $\mathbb{R}_{+}^n$ with $u\leq C$ and $|\nabla u|\leq C$ for some constant $C$ where $f$ is $C^1$ function and $\kappa$ is a non-zero constant. If there exists a nonpositive primitive $F$ of $f$ such that \cref{cond1} holds, then the solution $u$ is parallel.
\end{corollary}

\medskip

With respect to \Cref{prop2.1}, several remarks are due. The gradient estimate for the prescribed mean curvature equation has been studied extensively. The interior gradient estimate for the minimal surface equation was first obtained by Finn \cite{MR66533} (in dimension 2) and by Bombieri, De Giorgi and Miranda \cite{MR248647} (in all dimensions). By using the test function argument and a resulting Sobolev inequality, the interior gradient estimate for the general curvature equation had been obtained \cite{MR265745, MR324187, MR412605}. In 1983, Korevaar \cite{MR843597} introduced the normal variation technique for the interior gradient estimate on the minimal surface equation. In 1998, Wang \cite{MR1617971} gave another new proof for the estimate about the mean curvature equation by the Bernstein technique. The Dirichlet problem (in bounded domains) for the prescribed mean curvature equation has been studied by Jenkins and Serrin \cite{MR222467}. For the Neumann problem (in bounded domains) for the prescribed mean curvature equation (i.e., with constant contact angle condition on the boundary), Ural'tseva \cite{MR638359} first got the gradient estimates. Later, more general quasilinear equations in divergence form and more general boundary conditions were investigated \cite{MR708644, MR914813, MR3451945}. For unbounded domain, Colombo, Mari and Rigoli \cite{MR4271788} gave the gradient estimate for positive constant mean curvature graphs. Inspired by their method, we give the gradient estimate of capillary graphs whose mean curvature depends only on the height (i.e., solutions of problem \Cref{e.R}).

\medskip

{\bf Organization of the paper.} \Cref{sec2} will be devoted to the proof of \Cref{coro} by using a Modica type estimate. In Sections \ref{sec3} to \ref{secM} we will prove \Cref{Tmain}, and in order to simplify the notation we will suppose $\kappa = -1$. Observe that the constant $\kappa$ in \eqref{e.R} is clearly nonpositive, and for a general negative $\kappa$ the proof would be the same. In fact, by applying \Cref{coro} to a domain $\Omega$ satisfying the hypothesis of \Cref{Tmain}, we will get that $\Omega$ either is a half-plane or contains a half-plane. In this last case, we will be inspired by the strategy of \cite{MR3666566} to continue the proof. First, in \Cref{sec3} we will show the boundedness of the curvature of $\pp \Omega$ by using a blow-up strategy, that allows us to conclude thanks to a rigidity result for capillary minimal graphs by Cui \cite{MR4494613} (such result can be seen as the counterpart of the Traizet's result on harmonic overdetermined problem in the plane, \cite{MR3192039}.). Using the boundedness of the curvature of $\pp \Omega$, in \Cref{sec4} we will show that we can take the limit of our solution under translations obtaining a parallel solution in a half-plane. In \Cref{sec5} we obtain a characterization of the existence of a parallel solution in a half-plane with respect to the function $f$, and in \Cref{sec6} we will show that such parallel solution can be approached by a sequence of solution to \eqref{e.R} in large balls. This fact allows us to compare the parallel solution with our initial solution by using the Maximum Principle and obtain the conclusion in \Cref{secM}. The last part of the paper, \Cref{sec_boundedDU}, is devoted to the proofs of \Cref{prop2.1} and \Cref{cor_cap}.

\medskip


{\bf Acknowledgements.}
Y. L. and P. S. have been supported by the Grants PID2020-117868GB-I00 and PID2023-150727NB-I00 of the MICIN/AEI. P. S. has been supported also by the \emph{IMAG-Maria de Maeztu} Excellence Grant CEX2020-001105-M funded by the MICIN/AEI.


%
%

\section{Modica type estimate and proof of \Cref{coro}}\label{sec2}
This section is devoted to prove \Cref{coro}. We will obtain it as a corollary of a more general result. Let us consider a $C^1$ domain $\Omega \subset\mathbb{R}^{n}$, $n\geq 2$ (here $\Omega$ can be a bounded or an unbounded domain), and a solution $u$ of problem \eqref{e.R},
where $f$ is a given $C^1$ function and $\kappa$ is a constant. We will suppose that $u \in C^3(\Omega)$ and that $u\leq C$ and $|\nabla u|\leq C$ for some constant $C$.
Given a primitive $F$ of $f$, we define the $P$-function:
\begin{equation}\label{Pfunction}
P(x)=F(u(x))-(1+|\nabla u|^2)^{-\frac{1}{2}}.
\end{equation}
Based on the classical Modica estimate \cite{MR803255} and inspired by the Modica type estimate for elliptic overdetermined problems in the form \eqref{serrin} obtained by Ruiz, Wu and the second author in \cite{Ruiz_Sicbaldi_Wu24}, we will prove the following:

\begin{theorem} \label{Th.P1}
Let $F\in C^{2}(\mathbb{R})$ be a non-positive primitive of $f$ and $P$ be given by \cref{Pfunction}. Then
\begin{equation}\label{Modest}
  P(x)\leq\max\{-1,F(0)-(1+\kappa^2)^{-\frac{1}{2}}\} \ \mbox{ for all } x \in \Omega.
\end{equation}
Moreover, if there exists a point $x_{0}\in\Omega$ such that
\begin{equation*}
  P(x_{0})=\max\{-1,F(0)-(1+\kappa^2)^{-\frac{1}{2}}\},
\end{equation*}
then $P$ is constant, $\Omega$ is a half-space and $u$ is parallel.
\end{theorem}
\begin{remark}\label{re2.1}
In this theorem, the constant $\kappa$ can be $0$.
\end{remark}
\medskip

We will refer to estimate \eqref{Modest} as the {\it Modica type estimate for the general capillary overdetermined problems in the form \eqref{e.R}}. In the case where $P$ is bounded above by $F(0)-(1+\kappa^2)^{-\frac{1}{2}}$ we can give information on the mean curvature of $\partial \Omega$, as in the following statement:
\begin{proposition} \label{Th.P2}
Assume that $\kappa \neq 0$, and
\begin{equation} \label{cond}
P(x) \leq F(0)-(1+\kappa^2)^{-\frac{1}{2}} \ \mbox{ for all } x \in \Omega.
\end{equation}
Then, $H(p) \leq 0$ for any $p \in \partial\Omega$. Moreover, if there exists $p \in \partial \Omega$ such that $H(p)=0$, then $P$ is constant, $\Omega$ is either a half-space or a slab, that is, the domain between two parallel hyperplanes, and $u$ is parallel.
\end{proposition}

We just note that a slab is the part of the space between two parallel hyperplanes, i.e., $\mathbb{S} = \{x \in \R^n:  \ a \cdot (x-x_0) \in (z_1,z_2) \}$
for some $x_0,a \in \R^n$ and $z_1, z_2 \in \R$. A parallel solution in a slab is then a solution of the form
\[
u(x)=g(a\cdot (x-x_0)), \ \ x \in \mathbb{S}.
\]
for some function $g: (z_1,z_2) \to \R$.

\medskip

We will prove now the main result of this section, i.e., \Cref{Th.P1}. Later we will deduce \Cref{Th.P2} and from that \Cref{coro}. Suppose the assumptions of \Cref{Th.P1} are satisfied. We start with an important lemma, that shows that $P$ is a subsolution of an elliptic PDE. This will allow us to use later the maximum principle. In the following we will write $f$ instead of $f(u)$, as well as $f'$ instead of $f'(u)$. We use $x=(x_1,...,x_n)$ to denote a general vector in $\mathbb{R}^n$. In addition, $u_i$ denotes $\partial u/\partial x_i$ and $u_{ij}$ denotes $\partial^2u/\partial x_i\partial x_j$ etc. ($1\leq i,j\leq n$). We also use the common Einstein notation dropping the symbol of sum, that is understood when indices are repeated.

\begin{lemma}\label{Subsolution}
The function $P$ satisfies:
\begin{equation}\label{P-equation}
(1+|\nabla u|^2)\Delta P-u_iu_jP_{ij}+b^iP_i \geq 0,
\end{equation}
where
\begin{equation*}
  b^i=(1+|\nabla u|^2)^{\frac{1}{2}}(1+|\nabla u|^2+|\nabla u|^{-2})fu_i
+(1-|\nabla u|^{-2})u_ju_{ij},
\end{equation*}
for any $x\in \Omega$ such that $\nabla u(x)\neq 0$.
\end{lemma}
\proof Note that $P\in C^2$ since $u\in C^3$. By differentiating $P$, we have
\begin{equation}\label{P}
  \begin{aligned}
    &P_i=fu_i+(1+|\nabla u|^2)^{-\frac{3}{2}}u_ku_{ki}
  \end{aligned}
\end{equation}
and
\begin{equation}\label{P-2}
  \begin{aligned}
    &P_{ij}=f'u_iu_j+fu_{ij}+(1+|\nabla u|^2)^{-\frac{3}{2}} \left(-3(1+|\nabla u|^2)^{-1} u_ku_{ki}u_lu_{lj}
    +u_{ki}u_{kj}+u_ku_{kij} \right).
  \end{aligned}
\end{equation}
Then we have
\begin{equation}\label{P-form-1}
P_iu_i=f|\nabla u|^2+(1+|\nabla u|^2)^{-\frac{3}{2}}u_iu_ku_{ik},
\end{equation}
\begin{equation}\label{P-form-2}
P_iu_ju_{ij}=fu_iu_ju_{ij}+(1+|\nabla u|^2)^{-\frac{3}{2}}u_ku_{ki}u_ju_{ij},
\end{equation}
\begin{equation}\label{P-form-3}
|\nabla P|^2=f^2|\nabla u|^2+(1+|\nabla u|^2)^{-3}u_ku_{ki}u_ju_{ji}
  +2f(1+|\nabla u|^2)^{-\frac{3}{2}}u_iu_ku_{ki}
\end{equation}
and
\begin{equation}\label{P-form-4}
\Delta P=f'|\nabla u|^2+f\Delta u+(1+|\nabla u|^2)^{-\frac{3}{2}} \left( -3(1+|\nabla u|^2)^{-1} u_ku_{ki}u_lu_{li}
    +u_{ki}^2+u_k(\Delta u)_k \right).
\end{equation}
Hence, by \cref{P-form-1}, we get
\begin{equation}\label{u-3}
\begin{aligned}
  u_iu_ju_{ij}=&(1+|\nabla u|^2)^{\frac{3}{2}}\left(P_i u_i-f|\nabla u|^2\right).
\end{aligned}
\end{equation}
By \cref{P-form-2} and \cref{u-3}, we have
\begin{equation}\label{u-4}
\begin{aligned}
  u_ku_{ki}u_ju_{ji}=&(1+|\nabla u|^2)^{\frac{3}{2}}(P_iu_ju_{ij}-fu_iu_ju_{ij})\\
  =&(1+|\nabla u|^2)^{\frac{3}{2}}P_iu_ju_{ij}-(1+|\nabla u|^2)^{3}(fu_iP_i-|\nabla u|^2f^2).
\end{aligned}
\end{equation}
Inserting \cref{P-form-1} into \cref{P-form-3}, and combining with \cref{u-4}, we obtain
\begin{equation}\label{Pi_squre}
\begin{aligned}
  |\nabla P|^2=&2fu_iP_i-f^2|\nabla u|^2+(1+|\nabla u|^2)^{-3}u_ku_{ki}u_ju_{ji}\\
  =&fu_iP_i+(1+|\nabla u|^2)^{-\frac{3}{2}}P_iu_ju_{ij}.
  \end{aligned}
\end{equation}
By applying the H\"{o}lder inequality to \cref{P}, we get
\begin{equation*}
\begin{aligned}
|P_i-fu_i|=(1+|\nabla u|^2)^{-\frac{3}{2}}|u_ku_{ki}|
\leq (1+|\nabla u|^2)^{-\frac{3}{2}} |\nabla u| \left(\sum_{k=1}^{n} u_{ki}^2\right)^{\frac{1}{2}}.
\end{aligned}
\end{equation*}
Then
\[
u_{ki}u_{ki} \geq |\nabla u|^{-2} u_ku_{ki}u_ju_{ji}
\]
and using the first equality in \cref{Pi_squre}, we obtain
\begin{equation}\label{estimate-uki}
\begin{aligned}
u_{ki}u_{ki}\geq & (1+|\nabla u|^2)^3|\nabla u|^{-2}\left(|\nabla P|^2-2fu_iP_i+f^2|\nabla u|^2\right) \,.
\end{aligned}
\end{equation}
Write now the equation in \cref{e.R} in non-divergence form:
\begin{equation}\label{eq-div}
  (1+|\nabla u|^2)^{-\frac{1}{2}}\Delta u-(1+|\nabla u|^2)^{-\frac{3}{2}}u_iu_ju_{ij}+f(u)=0.
\end{equation}
Then, by differentiating with respect to $x_k$, we get
\begin{equation*}
\begin{aligned}
  -(1+|\nabla u|^2)^{-\frac{3}{2}}u_iu_{ik}\Delta u
  +(1+|\nabla u|^2)^{-\frac{1}{2}}(\Delta u)_k
  +3(1+|\nabla u|^2)^{-\frac{5}{2}}u_lu_{lk}u_i u_j u_{ij}\\
  -2(1+|\nabla u|^2)^{-\frac{3}{2}}u_{ik}u_ju_{ij}
  -(1+|\nabla u|^2)^{-\frac{3}{2}}u_iu_ju_{ijk}+f' u_k=0\,.
\end{aligned}
\end{equation*}
Next, multiply by $u_k$ both sides and take the sum over $k$:
\begin{equation}\label{diff-eq}
\begin{aligned}
  -(1+|\nabla u|^2)^{-\frac{3}{2}}u_i u_k u_{ik}\Delta u
  +(1+|\nabla u|^2)^{-\frac{1}{2}}(\Delta u)_k u_k
  +3(1+|\nabla u|^2)^{-\frac{5}{2}}u_l u_k u_{lk}u_i u_j u_{ij}\\
  -2(1+|\nabla u|^2)^{-\frac{3}{2}}u_ku_{ik}u_ju_{ij}
  -(1+|\nabla u|^2)^{-\frac{3}{2}}u_iu_ju_ku_{ijk}+f' |\nabla u|^2=0.
\end{aligned}
\end{equation}
Inserting \cref{P-form-1} into \cref{eq-div}, we obtain
\begin{equation}\label{Laplace-u}
\Delta u=(1+|\nabla u|^2)^{\frac{1}{2}}P_i u_i-(1+|\nabla u|^2)^{\frac{3}{2}}f.
\end{equation}
By considering \cref{P-form-4}$\times (1+|\nabla u|^2)-$\cref{P-2}$\times u_iu_j-$\cref{diff-eq}, we have
\begin{equation}\label{P-eq}
\begin{aligned}
(1+|\nabla u|^2)\Delta P-u_iu_jP_{ij}=&(1+|\nabla u|^2)f\Delta u-fu_iu_ju_{ij}
+(1+|\nabla u|^2)^{-\frac{3}{2}}u_i u_k u_{ik}\Delta u\\
&+(1+|\nabla u|^2)^{-\frac{1}{2}}u_{ki}^2-2(1+|\nabla u|^2)^{-\frac{3}{2}}u_ku_{ki}u_lu_{li}.
\end{aligned}
\end{equation}
There are five terms in the righthand. For the second and third terms, use \cref{eq-div} to replace $u_iu_ju_{ij}$; for the forth term, use \cref{estimate-uki} to replace $u_{ki}^2$; for the fifth term, use \cref{u-4} to replace $u_ku_{ki}u_lu_{li}$. Then
\begin{equation*}
  \begin{aligned}
    (1+|\nabla u|^2)\Delta P-u_iu_jP_{ij}\geq &f\Delta u-(1+|\nabla u|^2)^{\frac{3}{2}}f^2+
  (1+|\nabla u|^2)^{-\frac{1}{2}}(\Delta u)^2\\
  &+(1+|\nabla u|^2)^{\frac{5}{2}}|\nabla u|^{-2}\left(|\nabla P|^2-2fu_iP_i+f^2|\nabla u|^2\right)\\
  &-2P_iu_ju_{ij}+2(1+|\nabla u|^2)^{\frac{3}{2}}\left(fu_iP_i-f^2|\nabla u|^2\right).
  \end{aligned}
\end{equation*}
By inserting \cref{Pi_squre} into the above inequality, we have
\begin{equation*}
  \begin{aligned}
    (1+|\nabla u|^2)\Delta P-u_iu_jP_{ij}\geq &f\Delta u-(1+|\nabla u|^2)^{\frac{3}{2}}f^2+
  (1+|\nabla u|^2)^{-\frac{1}{2}}(\Delta u)^2\\
   +(1+|\nabla u|^2)^{\frac{3}{2}}
    &\left(|\nabla u|^{-2}-1\right)
  \left(|\nabla P|^2-2fu_iP_i+f^2|\nabla u|^2\right).
  \end{aligned}
\end{equation*}
Next, replace $\Delta u$ by \cref{Laplace-u} in the above inequality and after an arrangement, we obtain
\begin{equation*}
  \begin{aligned}
    (1+|\nabla u|^2)\Delta P-u_iu_jP_{ij}\geq &-(1+|\nabla u|^2)^{\frac{1}{2}}(1+2|\nabla u|^{-2})fu_iP_i\\
    &+(1+|\nabla u|^2)^{\frac{1}{2}}P_iu_iP_ju_j+(1+|\nabla u|^2)^\frac32(|\nabla u|^{-2}-1)|\nabla P|^2.
  \end{aligned}
\end{equation*}
Finally, by inserting \cref{Pi_squre} into the above inequality, we have
\begin{equation*}\label{P-eq2}
\begin{aligned}
&(1+|\nabla u|^2)\Delta P-u_iu_jP_{ij}+(1+|\nabla u|^2)^{\frac{1}{2}}(1+|\nabla u|^2+|\nabla u|^{-2})fu_iP_i
+(1-|\nabla u|^{-2})u_ju_{ij}P_i \\
&\geq (1+|\nabla u|^2)^{\frac{1}{2}}P_iu_iP_ju_j\geq 0.
\end{aligned}
\end{equation*}
This concludes the proof of the lemma.\proofend

Throughout the paper, we will denote by $B_\delta(p)$ the Euclidean ball of center $p \in \R^n$ and radius $\delta$. We will shortly use $B_\delta$ for $B_\delta(0)$, where $0$ is the origin of $\R^n$.
Next, we present the proof of the first part of \Cref{Th.P1}, i.e., the following:

\begin{proposition}
Under the assumptions of \Cref{Th.P1}, we have:
\begin{equation*}
 P(x) \leq \max\left\{-1,F(0)-(1+\kappa^2)^{-\frac{1}{2}}\right\},~~x\in \Omega.
\end{equation*}
\end{proposition}

\medskip

\proof Let
\begin{equation*}
  \bar{P}=\max\left\{-1,F(0)-(1+\kappa^2)^{-\frac{1}{2}}\right\}.
\end{equation*}
By the boundedness of $u$, $P$ is bounded. We prove this proposition by contradiction. Assume that $\bar{P}<\hat{P}:=\sup\limits_{x\in\Omega} P(x)$. Take a sequence of points $\{x_{k}\}_{k\in\mathbb{N}}\subset\Omega$ such that
\begin{equation*}
  P(x_{k})\rightarrow\hat{P} ~~\mbox{ as }~~ k\rightarrow \infty.
\end{equation*}

If $x_{k}$ is bounded, up to a subsequence we can assume that $x_k \to x_{0}$, with $P(x_0)=\hat{P}$. Since $\hat{P}>-1$ from the definition of $P$, $\nabla u(x_0)\neq 0$. Moreover, by $\hat{P} > F(0)-(1+\kappa^2)^{-\frac{1}{2}}$, we have $x_0 \in \Omega$. Hence, $P$ attains a local maximum at an interior point $x_0$. Since $P$ satisfies \cref{P-equation}, by the strong maximum principle, $P$ is constant in a neighborhood of $x_0$. This argument implies that the set
 $$ \{x \in \Omega: \ P(x)= \hat{P}\}$$
 is a non-empty open subset of $\Omega$. Since it is obviously closed, then $P\equiv \hat{P}$. But $P(x)=  F(0)-(1+\kappa^2)^{-\frac{1}{2}} < \hat{P}$ if $x \in \partial \Omega$, which is a contradiction.

We now prove the result in the case that $x_{k}$ is unbounded. We divide this part in two cases.

\textbf{Case 1:} $\limsup_{k \to +\infty} \mbox{dist}(x_{k},\partial\Omega)>\delta$ for some $\delta>0$.

Extend $u$ by $0$ outside $\Omega$ such that $u: \R^n \to \R$ is a globally Lipschitz function. Let $u_{k}(x)=u(x+x_{k})$ and $P_{k}(x)=P(x+x_{k})$. Then by the assumption,
\begin{equation}\label{eq12}
P_{k}(0)=P(x_{k})\to\hat{P}~~\mbox{as}~~ k\rightarrow \infty.
\end{equation}
Take a subsequence, still denoted by $u_k$, so that on compact sets of $\R^n$ the sequence $u_k$ converges uniformly to a certain Lipschitz function $u_\infty \geq 0$. Let
\begin{equation*}
\Omega_{\infty}:= \{x \in \R^n: \ u_\infty (x) >0\},~~P_{\infty}=F(u_{\infty})-(1+|\nabla u_{\infty}|^2)^{-\frac{1}{2}}.
\end{equation*}

Now, we show that $\Omega_{\infty}$ is non-empty. Indeed, since $\mbox{dist}(x_{k},\partial\Omega)>\delta$,
\begin{equation}\label{e2.1}
  \mathrm{div} \left(\frac{\nabla u_{k}}{\sqrt{1+|\nabla u_{k}|^2}}\right)+f(u_{k}) =0~~~~\mbox{ in}~~B_{\delta}.
\end{equation}
Since $\nabla u_k$ is uniformly bounded and $f$ is Lipschitz continuous, by the interior $C^{1,\alpha}$ estimate (see \cite[Theorem 13.1]{MR1814364}) and the $C^{2,\alpha}$ estimate (see  \cite[Theorem 6.2]{MR1814364}), $u_k$ is bounded in $C^{2,\alpha}$ norm  (for any $0<\alpha<1$) in any compact subset of $B_{\delta}$. Then $u_k$ converges to $u_{\infty}$ in $C^{2,\alpha}$ sense in compact subsets of $B_{\delta}$ and $u_{\infty}$ satisfies
\begin{equation}\label{e.case1.1}
  \mathrm{div} \left(\frac{\nabla u_{\infty}}{\sqrt{1+|\nabla u_{\infty}|^2}}\right)+f(u_{\infty}) =0~~~~\mbox{ in}~~B_{\delta}.
\end{equation}
In addition,
\begin{equation} \label{jo}
P_{\infty}(0)=\lim_{k\to \infty} P_k(0)=\hat{P}>\bar{P}.
\end{equation}
Since $\hat{P}>-1$, $\nabla u_{\infty}(0)\neq 0$. By combining with $u_{\infty}\geq 0$, we have $u_{\infty}(0)>0$ (otherwise $\nabla u_{\infty}(0)= 0$). Therefore, $0\in \Omega_{\infty}$ and then $\Omega_{\infty}$ is not empty.

For any $\Omega''\subset \subset \Omega'\subset\subset \Omega_{\infty}$, there exists $a>0$ such that $u_{\infty}\geq a$ in $\Omega'$. Since $u_k\to u_{\infty}$ uniformly, for $k$ large enough, $u_k>0$ in $\Omega'$. Then $u_k$ satisfies \cref{e2.1} in $\Omega'$. By a similar argument as above, we conclude that $u_k\to u_{\infty}$ in $C^{2,\alpha}(\bar{\Omega}'')$. Hence, $u_k\to u_{\infty}$ in $C^{2,\alpha}$ sense in any compact set of $\Omega_{\infty}$. Thus, $u_{\infty}$ satisfies \cref{e.case1.1} in $\Omega_{\infty}$ and $P_k\to P_{\infty}$ in $C^{1,\alpha}$ sense in any compact set of $\Omega_{\infty}$. Hence,
\begin{equation*}
P_{\infty}\leq \hat{P} ~~\mbox{ in}~~\Omega_{\infty}.
\end{equation*}

Note that $P_{\infty}(0)=\hat{P}$. By the strong maximum principle again, we conclude that
\begin{equation} \label{jo2}
P_{\infty} \equiv \hat{P} ~~\mbox{  in } \tilde{\Omega}_{\infty},
\end{equation}
where $\tilde{\Omega}_{\infty}$ is the connected component of $\Omega_{\infty}$ containing $0$.

Since $u_{\infty}$ is bounded in $\Omega_{\infty}$, we can choose a sequence of $y_k \in \tilde{\Omega}_{\infty}$ such that
\begin{enumerate}
	\item[i)] $u_{\infty}(y_k) \to \xi= \sup \{u_{\infty}(x):\ x \in \tilde{\Omega}_{\infty}\}>0$;
	\item[ii)]  $\nabla u_{\infty}(y_k) \to 0$ (observe that $y_k$ does not converge to $\partial \tilde{\Omega}_{\infty}$ since $u_\infty = 0$ on $\partial \tilde{\Omega}_{\infty}$).
\end{enumerate}
However, $\hat{P}=P_{\infty}(y_k) \to F(\xi)-1 \leq -1$, which is a contradiction.
\medskip

\textbf{Case 2:}  $\lim_{k \to +\infty} \mbox{dist}(x_{k},\partial\Omega)=0.$

As in Case 1, we extend $u$ by $0$ outside $\Omega$ such that $u: \R^n \to \R$ is a globally Lipschitz function. Since $u(x_{k})\rightarrow 0$, $P(x_k)\to \hat{P}$. Moreover, by the assumptions,
\begin{equation} \label{jo3}
\lim_{k\to\infty} |\nabla u(x_{k})|^{2} = \lim_{k\to\infty} \left(P(x_k)-F(u(x_k))\right)^{-2}-1 = \left(\hat{P} -F(0)\right)^{-2}-1 >\kappa^{2}.
\end{equation}
Denote $h_{k}:=\mbox{dist}(x_{k},\partial\Omega)\rightarrow0$ and assume that this distance is attained at $y_{k}\in\partial\Omega.$ Let
\[u_{k}(x):=\frac{1}{h_{k}}u(y_{k}+h_{k}x),~x\in \mathbb{R}^n.\]
Observe that $u_k(0)=0$ and $\nabla u_k$ is uniformly bounded. Then we can take a subsequence such that $u_k$ converges uniformly in compact sets of $\mathbb{R}^n$ to a limit function $u_\infty \geq 0$, which is Lipschitz continuous (but possibly unbounded).

Consider $\Omega_{\infty}= \{x \in \R^n: \ u_{\infty}(x) >0\}$ as above. Now, we show that $\Omega_{\infty}$ is non-empty. Let $z_{k}:=(x_{k}-y_{k})/h_{k}$. We can assume that $z_{k}$ converges to $z_{\infty}$ since $|z_{k}|=1.$ As in Case 1, $\nabla u_k$ is uniformly bounded and $u_k$ satisfies
\begin{equation*}
  \mathrm{div} \left(\frac{\nabla u_{k}}{\sqrt{1+|\nabla u_{k}|^2}}\right)+h_k \cdot f(h_k u_k)=0~~\mbox{ in }~~B_1(z_{k})\,.
\end{equation*}
With the aid of the interior estimates, we conclude that $u_k \to u_{\infty}$ in $C^{2,\alpha}$ sense in compact sets of $B_1(z_{\infty})$ and $u_{\infty}$ satisfies
\begin{equation}\label{e.case2.2}
  \mathrm{div} \left(\frac{\nabla u_{\infty}}{\sqrt{1+|\nabla u_{\infty}|^2}}\right)=0~~\mbox{ in }~~B_1(z_{\infty}).
\end{equation}
Note that $u_\infty \geq 0$ in $B_1(z_\infty)$ and by \cref{jo3},
\begin{equation}\label{eqg1}
a:=|\nabla u_{\infty}(z_\infty)|=\lim_{k\to\infty} |\nabla u_{k}(z_{k})|
=\lim_{k\to\infty} |\nabla u(x_{k})|> |\kappa|.
\end{equation}
By the strong maximum principle we conclude that $u_{\infty}>0$ in $B_1(z_\infty)$, that is, $B_1(z_\infty) \subset  \Omega_{\infty}$ and thus $\Omega_{\infty}$ is non-empty.

The same argument as above works similarly in $\Omega_{\infty}$ and then $u_k \to u_{\infty}$ in $C^{2,\alpha}$ sense in compact sets of $\Omega_{\infty}$. Hence $u_{\infty}$ satisfies \cref{e.case2.2} in $\Omega_{\infty}$.
By differentiating \eqref{e.case2.2} with respect to $x_l$, we have
\begin{equation}\label{e2.2}
  (a^{ij}(u_{\infty})_{lj})_i=0~~\mbox{ in}~~\Omega_{\infty},
\end{equation}
where
\begin{equation*}
  a^{ij}:=\frac{1}{\sqrt{1+|\nabla u_{\infty}|^2}}\left(\delta^{ij}-
  \frac{(u_{\infty})_{i}(u_{\infty})_{j}}{1+|\nabla u_{\infty}|^2}\right).
\end{equation*}
Let $v=|\nabla u_{\infty}|^2/2$. Multiply $(u_{\infty})_l$ in \cref{e2.2} and sum up from $l=1$ to $n$. Note that the eigenvalues of $a^{ij}$ are $1$ (multiplicity $n-1$) and $(1+|\nabla u_{\infty}|^2)^{-1}$ (see (10.7) in \cite{MR1814364}). Then
\begin{equation}\label{e2.3}
  (a^{ij}v_j)_i=a^{ij}(u_{\infty})_{li}(u_{\infty})_{lj}\geq 0.
\end{equation}
That is, $|\nabla u_{\infty}|^2$ is a subsolution of some uniformly elliptic equation.

Observe that
\begin{equation*}
  \begin{aligned}
F(0)-\left(1+|\nabla u_{\infty}(z_{\infty})|^{2}\right)^{-\frac{1}{2}}
=&\lim_{k\to \infty} F(h_ku_{k}(z_{k}))-\left(1+|\nabla u_{k}(z_{k})|^{2}\right)^{-\frac{1}{2}}\\
=&\lim_{k\to \infty} F(u(x_{k}))-\left(1+|\nabla u(x_{k})|^{2}\right)^{-\frac{1}{2}}\\
=&\lim_{k\to \infty} P(x_k)=\hat{P}.\\
  \end{aligned}
\end{equation*}
In addition, for all $ x \in B_1(z_\infty)$,
\begin{equation*}
  \begin{aligned}
F(0)-\left(1+|\nabla u_{\infty}(x)|^{2}\right)^{-\frac{1}{2}}
=&\lim_{k\to \infty} F(h_ku_{k}(x))-\left(1+|\nabla u_{k}(x)|^{2}\right)^{-\frac{1}{2}}\\
=&\lim_{k\to \infty} F(u(y_k+h_kx))-\left(1+|\nabla u(y_k+h_kx)|^{2}\right)^{-\frac{1}{2}}\\
=&\lim_{k\to \infty} P(y_k+h_kx)
\leq \hat{P}.\\
  \end{aligned}
\end{equation*}
As a consequence,
\begin{equation*}
|\nabla u_{\infty}(z_\infty)| \geq |\nabla u_{\infty} (x)|,~\forall ~ x\in B_1(z_\infty).
\end{equation*}
Since $|\nabla u_{\infty}|^2$ is a subsolution (see \cref{e2.3}), by the strong maximum principle
\begin{equation*}
|\nabla u_{\infty} |\equiv |\nabla u_{\infty}(z_\infty)|=a~~\mbox{  in}~~\tilde{\Omega}_{\infty},
\end{equation*}
where $a$ is given in \cref{eqg1} and $\tilde{\Omega}_{\infty}$ is the connected component of $\Omega_{\infty}$ containing $B_1(z_{\infty})$. From this and \cref{e.case2.2} we obtain that $u_{\infty}$ is harmonic in $\tilde{\Omega}_{\infty}$.

Now, we show that $u_{\infty}$ is a linear function. For any $x_0\in \tilde{\Omega}_{\infty}$, there exists a unit vector, say $e_n=(0,..,0,1)$, such that
\begin{equation*}
  \nabla u_{\infty}(x_0)\cdot e_n=|\nabla u_{\infty}(x_0)|=a.
\end{equation*}
Let $w=\frac{\partial u_{\infty}}{\partial x_n}$. Then $w$ is harmonic in $\tilde{\Omega}_{\infty}$ (since $u_{\infty}$ is harmonic) and
\begin{equation*}
w(x_0)=\frac{\partial u_{\infty}}{\partial x_n}(x_0)=\nabla u_{\infty}(x_0)\cdot e_n=a=|\nabla u_{\infty}(x)|\geq w(x),~\forall ~x\in \tilde{\Omega}_{\infty}.
\end{equation*}
By the strong maximum principle, $w\equiv a$. By combining with $|\nabla u_{\infty}|\equiv a$, we have
\begin{equation*}
\frac{\partial u_{\infty}}{\partial x_n}\equiv a,~~\frac{\partial u_{\infty}}{\partial x_i}\equiv 0,~~1\leq i\leq n-1.
\end{equation*}
Therefore, $u_{\infty}$ is linear in $\tilde{\Omega}_{\infty}$. Additionally, by noting $u_{\infty}(0)=0$, we have
\[u_{\infty}(x)=ax_{n}~~\mbox{ in}~~\tilde \Omega_{\infty}.\]
In particular, $\tilde \Omega_{\infty}$ is the upper half-space.

Next, we try to arrive at a contradiction. For any $\varepsilon\in(0,1)$, by the uniform convergence of $u_{k}$ to $u_{\infty}$,
\begin{equation*}
|u_{k}-ax_{n}|<\varepsilon^2~~\mbox{ in}~~B_1^+,
\end{equation*}
for all large enough $k$, where $B_1^+$ is the upper half-ball, and
\begin{equation*}
\frac{1}{\sqrt{1+|\nabla u_{k}|^2}}\left(\delta^{ij}-\frac{(u_{k})_i(u_k)_j}{1+|\nabla u_{k}|^2}\right)(u_k)_{ij}=  \mathrm{div} \left(\frac{\nabla u_{k}}{\sqrt{1+|\nabla u_{k}|^2}}\right)\leq \varepsilon~~\mbox{ in }~~\left\{x_n>\varepsilon\right\}\cap B_1^+.
\end{equation*}

\begin{figure}
  \centering
  \includegraphics[width=10cm]{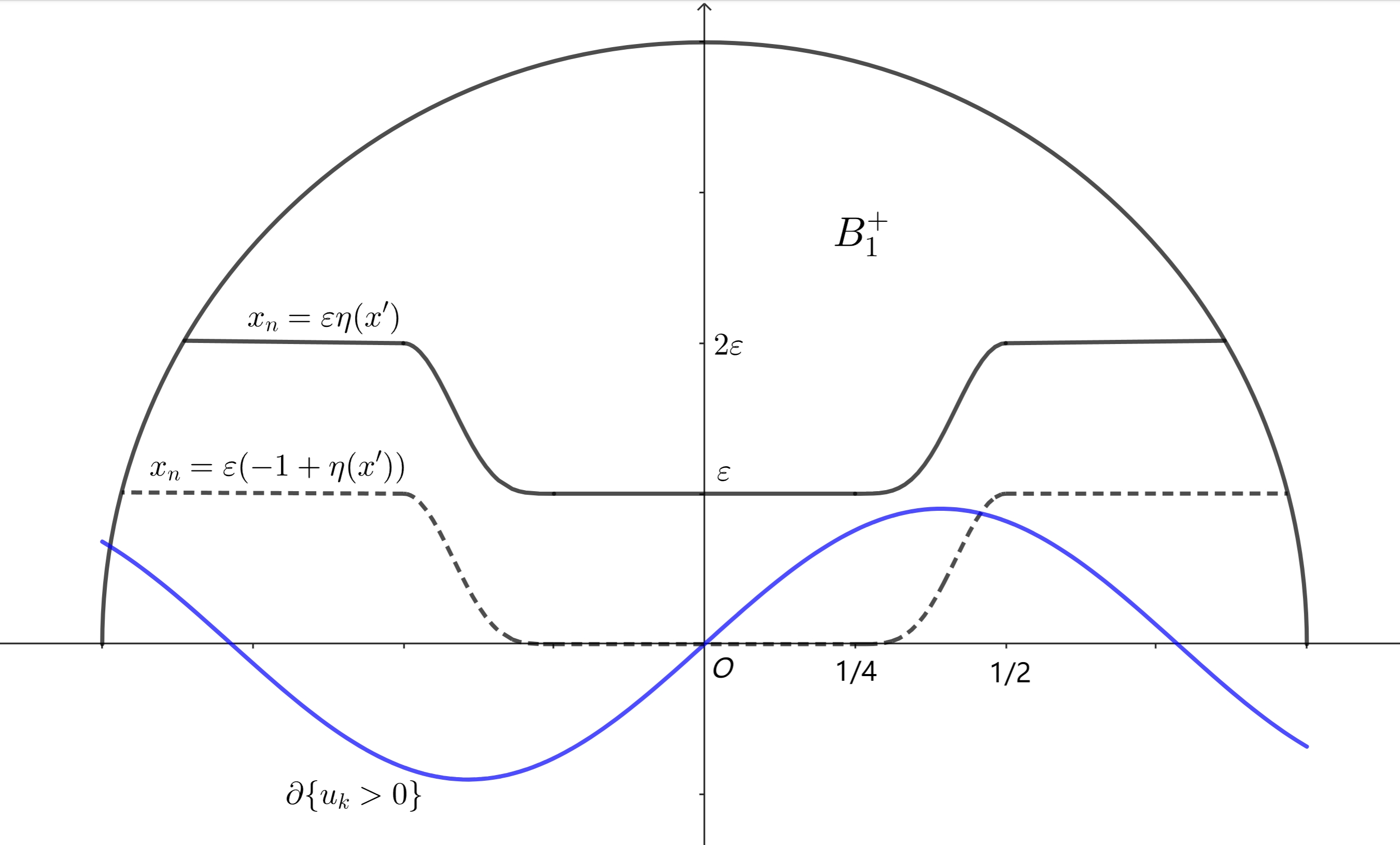}
  \caption{Schematic diagram of $D_t$}\label{PicDt}
\end{figure}

Let $\delta=(a-|\kappa|)/2$ and $1\leq \eta(x')\leq 2$ be a smooth function with
\begin{equation*}
  \begin{aligned}
    &\eta=1~~\mbox{on}~~|x'|\leq 1/4;\\
    &\eta=2~~\mbox{on}~~1/2\leq |x'|\leq1
  \end{aligned}
\end{equation*}
(here we denote $x = (x',x_n)$, where $x'=(x_1, ..., x_{n-1})$).
Consider the domain (see FIGURE \cref{PicDt})
\begin{equation*}
  D_{t}=\{x\in B_{1}^+:x_{n}>\varepsilon(t+\eta (x')),~~-1\leq t\leq 0\}
\end{equation*}
and set
\begin{equation*}
  v_t:=(a-\delta)\left(x_n-\varepsilon(\eta(x')+t)\right)+\frac{\delta}{2}\left(x_n-\varepsilon(\eta(x')+t)\right)^2.
\end{equation*}
By a direct calculation,
\begin{equation*}
  \begin{aligned}
    &(v_t)_{ij}=\varepsilon^2 \delta\, \eta_i\,\eta_j
    -\varepsilon\, \eta_{ij} \left((a-\delta)+\delta(x_n-\varepsilon(\eta+t))\right),~~~~i,j\leq n-1\,, \\
    &(v_t)_{in}=-\varepsilon \delta \eta_i,~~~~i\leq n-1\,,\\
    &(v_t)_{nn}=\delta.
  \end{aligned}
\end{equation*}
Note that $|\nabla u_k|$ is uniformly bounded. Hence, by choosing $\varepsilon$ small enough, we have for all large
$k$
\begin{equation*}
\frac{1}{\sqrt{1+|\nabla u_{k}|^2}}\left(\delta^{ij}-\frac{(u_{k})_i(u_k)_j}{1+|\nabla u_{k}|^2}\right)(v_t)_{ij}\geq \varepsilon.
\end{equation*}
In addition, if $\varepsilon$ small enough,
\begin{equation*}
  \begin{aligned}
    v_t
    &\leq a\left(x_n-\varepsilon(\eta(x')+t)\right)-\frac{\delta}{2}\left(x_n-\varepsilon(\eta(x')+t)\right)\\
    &=a x_n-\frac{\delta}{2} x_n-\varepsilon(\eta(x')+t)\left(a-\frac{\delta}{2}\right)\\
    &\leq a x_n-\frac{\delta}{2} x_n
    \leq a x_n- \varepsilon^2
  \end{aligned}
\end{equation*}
on $\partial B_1^+ \cap \partial D_{t}$.
In conclusion, we have for any $-1\leq t\leq 0,$
\begin{equation*}
 \begin{cases}
 \frac{1}{\sqrt{1+|\nabla u_{k}|^2}}\left(\delta^{ij}-\frac{(u_{k})_i(u_k)_j}{1+|\nabla u_{k}|^2}\right)(v_t)_{ij}\geq \varepsilon &\mbox{in $D_{t}$},\\
 v_{t}\leq u_k &\mbox{on $\partial B_1^+ \cap \partial D_{t}$},\\
 v_{t}=0&\mbox{on $ B_{1}^+\cap\{x_{n}=\varepsilon(t+\eta (x'))\}$}.
\end{cases}
\end{equation*}
Since $u_k(0)=0$, there exists $-1\leq t_0\leq 0$ and $p$ such that
\begin{equation*}
D_{t_{0}}\subset \{ u_k>0\}~~\mbox{ and }~~ p\in \{x_{n}=\varepsilon(t_0+\eta (x')),|x'|<1/2\}\cap \partial \{ u_k>0\}.
\end{equation*}
Thus, by the comparison principle, 
\begin{equation*}
|\kappa|=|(u_{k})_{\nu}(p)|\geq |(v_{t_0})_{\nu}(p)|
=|\nabla v_{t_0}(p)|\geq a-\frac{3}{2}\delta>|\kappa|,
\end{equation*}
where $\nu$ is the exterior unit normal at $p$. Hence, we arrive at a contradiction. \proofend

In order to complete the proof of \Cref{Th.P1}, we only need to show the following result. The proof idea is from \cite{MR1296785} (see also \cite{Ruiz_Sicbaldi_Wu24}).

\begin{proposition} \label{esta}
Under the assumptions of \Cref{Th.P1}, if $P(x)=\bar{P}$ at a point $x\in\Omega,$ then $P$ is constant, $u$ is parallel and $\Omega$ is a half-space.
\end{proposition}

\proof  First, we show that if $P(x_0)= \bar{P}$, then $\nabla u(x_0) \neq 0$. We prove it by contradiction. If $\nabla u(x_{0})=0,$ then $-1 \leq \bar{P}=F(u(x_{0}))-1 \leq -1$. Thus, $\bar{P}=-1$ and $F(u_0)=0$, where $u_0 = u(x_0)$. We consider the set
\[U=\{x\in \Omega:u(x)=u_{0}\}.\]
It is obvious that $U$ is closed and $U\neq\emptyset$. Let $x_{1}\in U$. For $t$ small, we consider the function $\varphi(t)=u(x_{1}+tw)-u_{0},$ where $w\in \mathbb{S}^{n-1}$ is arbitrarily fixed. Then we have
\[\varphi'(t)=\nabla u(x_{1}+tw)w.\]
Thus,
\begin{equation*}
      |\varphi'(t)|^{2}\leq|\nabla u(x_{1}+tw)|^{2}\leq \left(F(u(x_{1}+tw))-\bar{P}\right)^{-2}-1\,.
\end{equation*}
Define
\[
G(z) := \left(F(z)-\bar{P}\right)^{-2}-1\,.
\]
Note that $G$ attains the minimum at $u_0$ and $G(u_0)=\left(F(u_0)-\bar{P}\right)^{-2}-1=0$. Hence, $G'(u_0)=0$. Then we have $G(u)=O((u-u_{0})^{2})$ as $|u-u_{0}|\rightarrow0.$ Therefore,
   $$ |\varphi'(t)|\leq C|\varphi(t)|$$
for $t$ small enough. Since $\varphi(0)=0$, $\varphi\equiv0$ in $[-\varepsilon,\varepsilon]$ for some $\varepsilon>0$.  It shows that the set $U$ given above is open. Then $U=\Omega$. This implies that $u$ is a constant, which is a contradiction.

Now, we have that $\nabla u(x)\neq0$ for any $x\in\Omega$ with $P(x) = \bar{P}$. Next, we consider the following set
$$ \{ x \in \Omega: \ P(x)= \bar{P}\}\,.$$
By assumption it is non-empty. Applying the maximum principle to $P$, we conclude that the set is open. It is also closed by continuity. Then $P(x)= \bar{P}$ for all $x \in \Omega$. Hence, $\nabla P=0$. Then we obtain $\nabla P\cdot \nabla u=0$, i.e.,
\begin{equation*}
  F'(u)|\nabla u(x)|^{2}+\left(1+|\nabla u(x)|^{2}\right)^{-\frac{3}{2}}u_iu_ju_{ij}=0.
\end{equation*}
Combining this equality with the equation satisfied by $u$, we have
\begin{equation*}
  \Delta u=-\left(1+|\nabla u(x)|^{2}\right)^{\frac{3}{2}}F'.
\end{equation*}
Let $\tilde G\in C^{2}(\mathbb{R})$ and set $v=\tilde G(u)$. By the straightforward computation, we can obtain
\begin{equation} \label{otra}
   \begin{aligned}
   \Delta v=&\tilde G''(u)|\nabla u(x)|^{2}+\tilde G'(u)\Delta u\\
   =&\tilde G''(u)\left((F(u)-\bar{P})^{-2}-1\right)-\tilde G'(u)F'(u)\left(F(u)-\bar{P} \right)^{-3}.
   \end{aligned}
\end{equation}
Then we have that
\begin{equation}\label{eq41}
   \Delta v=0
\end{equation}
if we choose
   \[\tilde G(u)=\int_{u_{0}}^{u} \left((F(s)-\bar{P})^{-2}-1\right)^{-\frac{1}{2}}ds,\]
for some fixed $u_{0}\in u(\Omega).$ Observe that $(F(s)-\bar{P})^{-2}-1>0$ for any $u_0$ in the range of $u$, so that the integral defining $\tilde G$ is not singular. Moreover,
\begin{equation}\label{eq42}
    |\nabla v|^{2}=\tilde G'(u)^{2}|\nabla u|^{2}=1.
\end{equation}
Then we can infer that $v(x)=a\cdot x+b$ for some $a\in \mathbb{R}^{n}$ with $|a|=1$ and $b\in \mathbb{R}$. Notice that $\tilde G$ is invertible since it is increasing. Therefore, we can obtain that $u(x)=\tilde G^{-1}(v(x))=\tilde G^{-1}(a\cdot x+b)$. Then $u$ is parallel. A priori, $\Omega$ could be the inner space between two parallel hyperplanes, i.e., a slab, but this is impossible since $u$ has no critical points. This concludes the proof.
\proofend

Next, we prove \Cref{Th.P2} and \Cref{coro}.

\medskip

\noindent\textbf{Proof of \Cref{Th.P2}.}
Since $\kappa\neq 0$, we have that $u$ has no critical points close to $\partial \Omega$. By the assumption that $P$ attains its maximum at $\partial\Omega$, we have
    \[P_{\nu}\geq0\quad \mbox{on  $\partial\Omega.$}\]
Let $H$ be the mean curvature of $\partial\Omega$ at a given point. By \cref{Pfunction} and noting that $\nabla u=\kappa \nu$ on $\partial \Omega$, we have
\begin{equation}\label{e.P_nu}
    \begin{aligned}
      P_{\nu} &=f u_i \nu_i+\left(1+|\nabla u|^2\right)^{-\frac{3}{2}} u_{ij} u_i \nu_j\\
      &=\kappa f+\kappa \left(1+\kappa^2\right)^{-\frac{3}{2}}u_{ij}\nu_i \nu_j \geq 0
    \end{aligned}
\end{equation}
on $\partial \Omega$.
Moreover, the solution $u$ satisfies
\begin{equation}\label{e.equation.boundary}
  (1+\kappa^2)^{-\frac{1}{2}}\Delta u-\kappa^2 (1+\kappa^2)^{-\frac{3}{2}} u_{ij} \nu_i \nu_j+f=0
\end{equation}
on $\partial \Omega$.
In addition, by using the fact that $\frac{\partial^{2} u}{\partial\nu^{2}}=\Delta u-(n-1)H\frac{\partial u}{\partial\nu}$ on $\partial\Omega$ (see \cite[Section 5.4]{MR615561}), we have
\begin{equation}\label{e.H2}
  \Delta u- u_{ij}\nu_i \nu_j -\kappa (n-1)H=0.
\end{equation}
Then, by considering \cref{e.P_nu} $-\kappa \times$ \cref{e.equation.boundary}
$+\kappa (1+\kappa^2)^{-\frac{1}{2}} \times $\cref{e.H2}, we have
  \[H\leq0\]
on $\partial \Omega$ and that if $H(p)=0$ for some $p \in \partial \Omega$, then $P_{\nu}(p)=0$. In this last case, by Hopf's lemma, we can get that
$P$ is constant, at least in a neighborhood of $p$. We can now argue as in the proof of \Cref{esta} (see \cref{otra}, \cref{eq41}, \cref{eq42}) to conclude that $u$ is parallel in such neighborhood. By unique continuation, $u$ is parallel and $\Omega$ is either a half-space or a slab. ~\qed~\\

\noindent\textbf{Proof of \Cref{coro}.} By \Cref{Th.P1}, $P(x) \leq F(0)-\left(1+\kappa^{2}\right)^{-\frac{1}{2}}$ for all $x \in \Omega$.

Since $\kappa \neq 0$, according to \Cref{Th.P2}, $H \leq 0$. Moreover, if $H(p)=0$ at some point of $\partial \Omega$, then $P$ is constant and $u$ is parallel. Additionally, by \Cref{esta}, we have that $u$ has no critical points, and then $\Omega$ is a half-space.
~\qed~\\

\section{Boundedness of the curvature}\label{sec3}
In Sections \ref{sec3} to \ref{secM}, in order to simplify the notation, we will suppose $\kappa = -1$ (for a general $\kappa \neq 0$ the proof would be the same).

Let $\Omega \subset \R^2$ and $u$ satisfying the hypothesis of \Cref{Tmain}. In this section we will establish the boundedness of the curvature of the boundary of the domain. Before showing the result of the boundedness of the curvature, we recall a regularity result. In fact, in our assumptions the boundary of the domain is of class $C^{1}$. A regularity argument for elliptic equations shows that the solution of \Cref{e.R} is $C^{2,\alpha}$ and $C^{2,\alpha}$ up to the boundary. That is, the following regularity result holds (see \cite[Theorem 1]{MR1200301}):

\begin{lemma}\label{Le.Regularity}
Let $p \in \partial \Omega$, $\varepsilon>0$ such that $\Omega \cap B_\varepsilon(p)$ and $\partial \Omega \cap B_\varepsilon(p)$ are both connected. Then $u\in C^{2, \alpha}(\bar\Omega\cap B_r)$ and $\partial \Omega\cap B_r\in C^{2,\alpha}$ for any $0<r<\varepsilon$ and any $0 < \alpha < 1$. Moreover
\begin{equation*}
  \|u\|_{C^{2,\alpha}(\bar\Omega\cap B_{r})}\leq C,\quad \|\partial \Omega\cap B_r\|_{C^{2,\alpha}}\leq C,
\end{equation*}
where $C$ depends only on $\alpha,r,\varepsilon$, $\|u\|_{C^{0,1}(\bar\Omega \cap B_\varepsilon(p))}$, $\|f (u)\|_{C^{\alpha}}$ and $\|\partial \Omega\cap B_\varepsilon(p)\|_{C^{1}}$.
\end{lemma}

The previous lemma guarantees that our domain $\Omega$, originally supposed to be of class $C^1$, is in fact more regular, i.e., locally of class $C^{2,\alpha}$, and the solution $u$, originally supposed of class $C^1$ up to the boundary, is in fact locally of class $C^{2,\alpha}$ up to the boundary.

\medskip

The main result of this section is as follows.
\begin{proposition}\label{curv}
Let $\Omega \subset \R^2$ and $u$ satisfying the hypothesis of \Cref{Tmain}. Then:
\begin{enumerate}
\item[i)] The curvature of $\partial \Omega$ is bounded.
\item[ii)] The $C^{2,\alpha}$ norm of the function $u$ is bounded
in $\overline{\Omega}$ for any $\alpha>0$.
\end{enumerate}
\end{proposition}

\proof If i) holds, using \cite[Lemma 2.5]{MR3666566}, we have that in \Cref{Le.Regularity} the value $\varepsilon$ can be chosen uniformly for any $p \in \partial \Omega$, and so we have ii).

\medskip

We now turn our attention to the proof of i). We prove it by contradiction. Suppose that $\partial \Omega$ has unbounded curvature. First note that by \Cref{coro}, $(\bar\Omega)^c$ is a unbounded convex domain. Let $K(q)$ denote the curvature of $\partial \Omega$ at the point $q \in \partial \Omega$. If $K$ is unbounded, there exists a sequence of points $q_n \in \partial \Omega$ such that $|q_n|,|K(q_n)|\rightarrow +\infty$.


Let $I_n$ be the connected component of $\partial \Omega \cap B_1(q_n)$ containing $q_n$ and let $p_n\in I_n$ be the point where the function
\begin{equation*}
p \to d(p, \partial B_1(q_n)) \, |K(p)| = (1-|p-q_n|)\,
|K(p)|, \quad p \in I_n
\end{equation*}
attains its maximum, that clearly exists. Then we have
\begin{equation*}
 |K(p)|\leq \frac{1-|p_n-q_n|}{1-|p-q_n|}|K(p_n)|, \quad p \in I_n.
\end{equation*}

Let $r_n=1-|p_n-q_n|$ and $R_n=|K(p_n)|(1-|p_n-q_n|)=r_n|K(p_n)|$. Clearly $R_n\geq |K(q_n)|\to \infty$. Consider the transformation $T_n$ in $\R^2$ given by
\begin{equation*}
  z=|K(p_n)|\cdot(p-p_n)
\end{equation*}
where $z=(z_1,z_2)$ and a rotation such that $\Omega_n=T_n(\Omega)$ tangents to the $z_1$-axis at $0$. Next, we set $v_n(z)=|K(p_n)|\cdot u(p)$ and then $v_n$ satisfies
\begin{equation}\label{e.3.1.vn}
  \begin{cases}
  \mathrm{div} \left(\frac{\nabla v_n}{\sqrt{1+|\nabla v_n|^2}}\right)+f_n(v_n) =0& \mbox{in $\Omega_n\cap B_{R_n}$, }\\
  v_n>0 &\mbox{in $\Omega_n \cap B_{R_n}$, }\\
  v_n=0 &\mbox{on $\partial \Omega_n \cap B_{R_n}$, }\\
  \partial_{\nu} v_n=-1 &\mbox{on $\partial \Omega_n \cap B_{R_n}$, }
  \end{cases}
\end{equation}
where
\begin{equation*}
f_n (v_n)=\frac{f(u)}{|K (p_n)|}=\frac{1}{|K (p_n)|} f\left(\frac{v_n}{|K (p_n)|}\right),\quad \partial \Omega_n=|K (p_n)|\cdot \partial \Omega,\quad K_n (z)=\frac{K (p)}{|K (p_n)|},
\end{equation*}
being $K_n$ the curvature of $\partial \Omega_n$.
Then, we have $|K_n (0)|=|K(p_n)|/|K(p_n)|=1$ and
\begin{equation*}
  |K_n(z)|\leq \frac{r_n}{r_n-|z|/|K(p_n)|}|K_n(0)|
  \leq \frac{1}{1-|z|/R_n},\quad z\in \partial \Omega_n \cap B_{R_n}.
\end{equation*}
Moreover,
\begin{equation*}
  |\nabla v_n|=|\nabla u|\leq C.
\end{equation*}

For any fixed $R>0$, if $n$ is large enough we have $R\leq R_n/2$, and then
\begin{equation}\label{e3.1}
  |K_n(z)|
  \leq \frac{1}{1-R/R_n}
  \leq 2,\quad z\in \partial \Omega_n \cap B_R.
\end{equation}
Since $\|v_n\|_{L^{\infty}(\Omega_n \cap B_{R})}\leq CR$ because its gradient is bounded, by \Cref{Le.Regularity} we have
\begin{equation*}
  \begin{aligned}
    &\|v_n\|_{C^{2,\alpha}(\Omega_n \cap B_{R/2})}\leq
    C\left(\|v_n\|_{L^{\infty}(\Omega_n \cap B_{R})}+\|f_n\|_{C^{0,1}(\Omega_n \cap B_{R})}\right)
    \leq C,\\
  \end{aligned}
\end{equation*}
where $C$ depends on $R$ and is independent of $n$. That is, $v_n$ is bounded in the $C^{2,\alpha}$ sense on compact sets. Then $v_n$ converges, up to a subsequence, in $C^{2,\alpha}$ sense on compact sets to a solution $v_{\infty}$ of the problem
\begin{equation*}
  \begin{cases}
  \mathrm{div} \left(\frac{\nabla v_{\infty}}{\sqrt{1+|\nabla v_{\infty}|^2}}\right) =0& \mbox{in $\Omega_{\infty}$, }\\
  v_{\infty}>0 &\mbox{in $\Omega_{\infty}$, }\\
  v_{\infty}=0 &\mbox{on $\partial \Omega_{\infty}$, }\\
  \partial_{\nu} v_{\infty}=-1 &\mbox{on $\partial \Omega_{\infty}$\,, }
  \end{cases}
\end{equation*}
and $0 \in \partial \Omega_{\infty}$.
By the rigidity of this overdetermined problem about the minimal surface equation in the plane (see \cite[Theorem 1.1]{MR4494613}), we conclude that $\Omega_{\infty}$ is a half-plane. If $K_\infty$ denotes the curvature of $\partial \Omega_{\infty}$, we have $|K_\infty(0)|=0$ because $\Omega_{\infty}$ is a half-plane. On the other hand, $|K_\infty(0)| = \lim_{n \to +\infty} |K_n(0)| = 1$, which gives a contradiction. \proofend

\begin{remark}\label{re.3.1}
Comparing the previous result with \cite{MR3666566}, it is worth to point out that the same blow up process can by applied to problem \cref{bcnn}, and in such case one can obtain the uniform $C^{2,\alpha}$ bound for $v_n$ by the regularity theory for the Poisson's equation since the curvature of $\partial \Omega_n$ is bounded (see \cite[(3.1)]{MR3666566} and the proof following it). For our problem \cref{e.R}, without the assumption $|\nabla u|\leq C$ one cannot get the uniform $C^{2,\alpha}$ bound for $v_n$ even if the curvature of $\partial \Omega_n$ is bounded.
\end{remark}

\section{Building a parallel solution from the given overdetermined problem}\label{sec4}
Let $\Omega \subset \R^2$ and $u$ satisfying the hypothesis of \Cref{Tmain}. As in the previous section, we assume $\kappa = -1$. In this section, we build a parallel solution of \eqref{e.R} in a half-plane. Let $\partial \Omega$ be parameterized by $\gamma (t)$, where $t\in \mathbb{R}$ is the arc length parameter. After a rigid motion in the plane we can assume that the orthonormal basis $\{\nu(t), \gamma'(t)\}$ given by the exterior unit normal and the tangent vector is positive.

By \Cref{coro}, $(\bar{\Omega})^c$ is a unbounded convex domain. Then the limit of $\gamma'(t)$ exists as $t\to +\infty$. Without loss of generality, we assume
\begin{equation*}
\gamma(0)=(0,0),~~ \gamma'(0) \cdot (0,1) \geq 0, ~~\lim_{t\to +\infty} \gamma'(t)=(1,0).
\end{equation*}

Let $x,y$ the two variables of $\R^2$. We can show the:
\begin{proposition}\label{parallel-u}
There exist sequences $t_n \to +\infty$ such that, if $q_n=\gamma(t_n)$, we have
\begin{equation}\label{e4.1}
\left\{\begin{aligned}
&\Omega_n := \Omega-q_n\to \mathbb{R}^2_+~~\mbox{ locally in the Hausdorff distance},\\
&u_n(x,y):=u((x,y)-q_n)\to u_{\infty}\mbox{ locally in } \mathbb{R}^2_+,
\end{aligned}\right.
\end{equation}
where $u_{\infty}$ is a bounded parallel solution of \cref{e.R} in $\mathbb{R}^2_+$. Precisely,
\begin{equation*}
\left\{\begin{aligned}
&d_{H}(\Omega_n\cap B_R,B_R^+) \to 0,~\forall ~R>0,\\
&u_n\to u_{\infty}~~\mbox{ in }~~ C^{2}(\bar\Omega'),~\forall ~\Omega'\subset\subset \mathbb{R}^2_+,
\end{aligned}\right.
\end{equation*}
where $d_H$ denotes the Hausdorff distance.
\end{proposition}

\proof Since $\gamma'(t)\to (1,0)$ as $t\to +\infty$, for any $\varepsilon>0$, there exists $t_0>0$ such that if $t>t_0$,
\begin{equation}\label{vr_q}
  |\gamma'(t)-(1,0)|<\varepsilon/2.
\end{equation}

Take $t_1=t_0+8/\varepsilon$ and $q=(x^q,y^q)=\gamma(t_1)\in \partial \Omega$. Define the sector
\begin{equation*}
  G_{q,\varepsilon}=\left\{(x,y):\  y \leq y^q -\varepsilon |x-x^q| \right\}.
\end{equation*}
Let $I_{q,\varepsilon}$ be the connected component of $\partial \Omega\cap B_{1/\varepsilon}(q)$ containing $q$.
We claim that
\begin{equation}\label{e4.2}
I_{q,\varepsilon}\cap \partial G_{q,\varepsilon}=\{q\}.
\end{equation}
Suppose not, i.e., there exists another point $q_2=\gamma(t_2)\in I_{q,\varepsilon}\cap \partial G_{q,\varepsilon}$. Since $(\bar{\Omega})^c$ is convex, the length of $I_{q,\varepsilon}$ is less than the perimeter of $ B_{1/\varepsilon}(q)$, i.e. $2\pi/\varepsilon$. From the choice of $q$, we have $t_2>t_0$ since $2\pi/\varepsilon < 8/\varepsilon$. Then by the Lagrange mean value theorem, there exists $t\in (t_2,t_1)$ (or $t\in (t_1,t_2)$) such that $\gamma'(t)=(\cos \varepsilon,\sin \varepsilon)$ (or $\gamma'(t)=(\cos \varepsilon,-\sin \varepsilon)$), which contradicts \cref{vr_q} if we take $\varepsilon$ small enough.

From the above argument, we can choose sequences $\varepsilon_n\to 0$ and $t_n\to +\infty$ such that \cref{e4.2} holds for $\varepsilon_n$ and $q_n=\gamma(t_n)$. Let
\begin{equation*}
\Omega_n= \Omega-q_n,\quad u_n(x,y)=u((x,y)-q_n).
\end{equation*}
By \Cref{curv}, the curvature of $\Omega_n$ and the $C^{2,\alpha}$ norm of $u_n$ in $\overline{\Omega}_n$ are uniformly bounded. Then the sequence $\Omega_n$ converges in the Hausdorff distance to a domain $\Omega_{\infty}$ with boundary unbounded and connected and the sequence $u_n$ converges in $C^{2}$ sense to a solution $u_{\infty}$ defined in $\Omega_\infty$. Clearly $(0,0) \in \partial \Omega_\infty$ and, by \cref{e4.2}, $\Omega_{\infty}$ is contained in the upper half-plane. Then the curvature of $\partial \Omega_{\infty}$ at $(0,0)$ is nonnegative. With the aid of \Cref{coro}, $\Omega_{\infty}$ must be the half-plane, and then the solution $u_{\infty}$ is a bounded parallel solution. \proofend

\section{Characterization of bounded positive parallel solutions}\label{sec5}
In this section, we will show that the existence of a parallel solution to \eqref{e.R} with $\kappa = -1$ in a half-space, bounded and with bounded gradient, is equivalent to some properties on $f$ (given in term of its primitive $F$). First, we show the following:
\begin{lemma}\label{th1.2}
Let $\varphi$ be a solution of
\begin{equation}\label{eq.one}
  \left\{
  \begin{aligned}
&\frac{\varphi''(t)}{\left(1+|\varphi'(t)|^2\right)^{\frac{3}{2}}}+f(\varphi(t))=0~~\mbox{ in }~\mathbb{R}^+;\\
&\varphi(0)=0,~ \varphi'(0)=1.
\end{aligned}
  \right.
\end{equation}
Assume that there exist $L>0$ and a primitive $F$ of $f$ such that
\begin{equation*}
f(L)=0,~~F(L)=0> F(l),~\forall ~l<L~~\mbox{ and }~~F(0)=\frac{\sqrt{2}}{2}-1.
\end{equation*}
Then
\begin{equation*}
\begin{aligned}
&\varphi'>0~~\mbox{ in}~\mathbb{R}^+,~
\lim_{t\to +\infty} \varphi'(t)=\lim_{t\to +\infty} \varphi''(t)=0~~\mbox{ and }~~
\lim_{t\to +\infty} \varphi(t)=L.
\end{aligned}
\end{equation*}
\end{lemma}
\proof Define the Hamiltonian
\begin{equation}\label{e.H}
\mathcal{H}(t):=F(\varphi(t))-\frac{1}{\sqrt{1+|\varphi'(t)|^2}}.
\end{equation}
By the equation in \cref{eq.one}, $\mathcal{H}'\equiv 0$. Note that $\mathcal{H}(0)=-1$ and hence $\mathcal{H}\equiv -1$.

\textbf{Step 1.} $\varphi(t)<L$ for any $t>0$. Otherwise, there exists $t_0>0$ such that
$\varphi(t_0)=L$. By \cref{e.H} and $\mathcal{H}\equiv -1$, $\varphi'(t_0)=0$. Then by noting $f(L)=0$, we have that $\varphi\equiv L$ is the unique solution of
\begin{equation*}
  \left\{
  \begin{aligned}
&\frac{\varphi''(t)}{\left(1+|\varphi'(t)|^2\right)^{\frac{3}{2}}}+f(\varphi(t))=0~~\mbox{ in }~\mathbb{R}^+;\\
&\varphi(t_0)=L,~ \varphi'(t_0)=0,
\end{aligned}
  \right.
\end{equation*}
which is impossible.

\textbf{Step 2.} $\varphi'(t)\neq 0$ for any $t>0$. Suppose not, i.e, there exists $t_0>0$ such that
$\varphi'(t_0)=0$. By Step 1, $\varphi(t_0)<L$. Then
\begin{equation*}
0>F(\varphi(t_0))=-1+\frac{1}{\sqrt{1+|\varphi'(t_0)|^2}}=0,
\end{equation*}
which is a contradiction.

\textbf{Step 3.} $\varphi'(t)>0$ for any $t>0$. This is an immediate consequence of $\varphi'(0)=1$ and Step 2.

\textbf{Step 4.} $\lim_{t\to +\infty} \varphi(t)=L$. By Step 3, $\varphi$ is strictly increasing. Note that $\varphi(t)<L$ for any $t>0$. Hence, there exists $\tilde{L}$ such that
\begin{equation*}
  \lim_{t\to +\infty} \varphi(t)=\tilde{L}\leq L.
\end{equation*}
In addition,
\begin{equation*}
  \lim_{t\to +\infty} \varphi'(t)=  \lim_{t\to +\infty} \sqrt{\frac{1}{\left(F(\varphi(t))+1\right)^2}-1}
  =\sqrt{\frac{1}{\left(F(\tilde L)+1\right)^2}-1}.
\end{equation*}
Note that this limit must be 0. Otherwise, it would contradict the fact that $\varphi$ is bounded. Hence, $\tilde{L}=L$.

Finally, note that
\begin{equation*}
  \lim_{t\to +\infty}\varphi''(t)=-\lim_{t\to +\infty}f(\varphi(t))\left(1+|\varphi'(t)|^2\right)^{\frac{3}{2}}=-f(L)=0.
\end{equation*}
This concludes the proof of the result.
\proofend

Next, we prove the inverse.
\begin{lemma}\label{th1.4}
Let $\varphi$ be a solution of \cref{eq.one}. Assume that
\begin{equation}\label{e1.2.1}
\varphi~~\mbox{is positive and bounded in}~\mathbb{R}^+.
\end{equation}
Then there exist $L>0$ and a primitive $F$ of $f$ such that
\begin{equation}\label{e1.2.2}
f(L)=0,~F(L)=0> F(l),~\forall ~l<L,~F(0)=\frac{\sqrt{2}}{2}-1~~\mbox{ and }~~\inf_{l\leq L} F(l)>-1.
\end{equation}
\end{lemma}

\proof Since $\varphi >0$ in $\R^+$, by the moving plane (line) method, we have $\varphi'>0$ in $\R^+$ immediately. Thus, $\varphi$ is monotone increasing and bounded. Then there exists $L>0$ such that $\varphi(t)\rightarrow L$ as $t\rightarrow +\infty$.

Choose $F$ such that $F(0)=\sqrt{2}/2-1$. Note that
\begin{equation}\label{e5.1}
  \mathcal{H}(t)=F(\varphi(t))-\frac{1}{\sqrt{1+|\varphi'(t)|^2}}\equiv -1.
\end{equation}
Hence,
\begin{equation*}
  \lim_{t\to +\infty} \varphi'(t)=  \lim_{t\to +\infty} \sqrt{\frac{1}{\left(F(\varphi(t))+1\right)^2}-1}
  =\sqrt{\frac{1}{\left(F(L)+1\right)^2}-1}.
\end{equation*}
The limit must be $0$. Otherwise it would contradict the fact that $\varphi$ is bounded. Hence, $F(L)=0$. Furthermore,
\begin{equation*}
\lim_{t\to \infty} \varphi''(t)
=-\lim_{t\to \infty}f(\varphi(t))\left(1+|\varphi'(t)|^2\right)^{\frac{3}{2}}=-f(L).
\end{equation*}
Similarly, the limit must be $0$. That is, $f(L)=0$.

From \cref{e5.1} and $\varphi'$ being bounded (since $\varphi'(0)=1$ and $\varphi'\to 0$), we have $\inf F>-1$. Finally, for any $l\in (0,L)$, since $\varphi'>0$, there exists a unique $t_0>0$ such that $\varphi(t_0)=l$. By \cref{e5.1} and $\varphi'(t_0)>0$, we have $F(l)<0=F(L)$.\proofend


Let $\Omega \subset \R^2$ and $u$ satisfying the hypothesis of \Cref{Tmain}. From \Cref{parallel-u} and \Cref{th1.4} we have that there exists $L>0$ such that $f(L)=0$ and up to change the primitive $F$ we have:
\begin{equation*}
F(L)=0 > F(l),~\forall ~l<L~~,~~F(0)=\frac{\sqrt{2}}{2}-1~~\mbox{ and }~~\inf_{l\leq L} F(l)>-1.
\end{equation*}

\section{Existence of solutions in balls and asymptotic properties}\label{sec6}
In this section, we prove that if there exists a parallel positive solution of \eqref{e.R} in a half-space of $\R^n$, i.e., if there exists a bounded positive solution of the ODE \eqref{eq.one}, then there exist positive radial solutions of \eqref{e.R} in balls for any radius large enough. This result is valid in $\R^n$ and then will be stated in its general form.

Suppose that there exists a bounded positive solution of \eqref{eq.one}, where $f$ is a $C^1$ function.
Then by \Cref{th1.4}, there exist $L>0$ and $F$ (one of the primitive of $f$) such that
\begin{equation}\label{e6.1}
f(L)=0,~F(L)=0> F(l),~\forall ~l<L,~F(0)=\frac{\sqrt{2}}{2}-1~~\mbox{ and }~~\inf_{l\leq L} F(l)>-1.
\end{equation}
In fact, $F$ can be written as
\begin{equation*}
F(u)=\int_{0}^{u} f(s)\, ds-\int_{0}^{L} f(s)\, ds=\int_{0}^{u} f(s)\, ds+\frac{\sqrt{2}}{2}-1.
\end{equation*}
From now on, $L$ and $F$ are fixed as above.

The main result in this section is the following:
\begin{proposition}\label{th6.1}
Suppose that there exists a bounded positive solution $\varphi$ of \eqref{eq.one}, where $f$ is a $C^1$ function. Then
there exists $R_0>0$ such that for any $R>R_0$, there exists a positive radial solution $u_{R}\in C^{2,\alpha}(\bar{B}_R)$ (for any $0<\alpha<1$) of
\begin{equation*}
\begin{cases}
\mathrm{div}\left(\frac{\nabla u_{R}}{\sqrt{1+|\nabla u_{R}|^2}}\right)+f(u_{R})=0\quad~~&\mbox{in}~~B_R,\\
u_{R}=0\quad~~&\mbox{on}~~\partial B_R,
\end{cases}
\end{equation*}
(where $B_R$ is a ball of radius $R$ in $\R^n$).
In addition, we have the following estimates:
\begin{equation*}
0<u_{R}< L~~\mbox{ in}~~ B_R,\quad |\nabla u_{R}|^2\leq \frac{1}{(\inf F+1)^2}-1 \quad\mbox{in}~~B_R,\quad \|u_{R}\|_{C^{2,\alpha}(\bar{B}_R)}\leq C,
\end{equation*}
where $C$ depends only on $n, \alpha, L,\inf F$ and $\|f\|_{C^{0,1}}$.
\end{proposition}

In the following proof, we need the Pohozaev's identity (see \cite[Proposition 2]{MR2921887}), that we recall in the following:
\begin{lemma}\label{le1.1}
Let $\Omega\subset \mathbb{R}^n$ be a smooth bounded domain and $u\in C^2(\Omega)\cap C(\bar{\Omega})$ be a solution of
\begin{equation*}
  \begin{cases}
    \mathrm{div}\left(\varphi'(|\nabla u|^2)\nabla u\right)+g(u)=0~~~~&\mbox{ in }~~\Omega;\\
    u=0~~~~&\mbox{ on}~~\partial \Omega,
  \end{cases}
\end{equation*}
where $g \in C^1(\R)$ and $\varphi \in C^2(\R_+)$ is such that $\varphi(0)=0$ and $\varphi'(s) \geq 0$ for any $s \geq 0$. Suppose that $G$ is a primitive of $g$ with $G(0)=0$ and let $x_0\in \Omega$. Then
\begin{equation}\label{e1.2}
\frac{n}{2}\int_{\Omega} \psi(|\nabla u|^2) -n \int_{\Omega}G(u)
+\frac{1}{2} \int_{\partial\Omega}(x-x_0)\cdot \nu\, \tilde{\psi} (u_{\nu}^2)=0,
\end{equation}
where
\begin{equation*}
\psi(s)=\varphi(s)-\frac{2}{n} s\varphi'(s),\quad\tilde{\psi}(s)=2 s\varphi'(s)-\varphi(s)
\end{equation*}
(and $\nu$ is as usual the exterior normal vector on $\partial \Omega$).
\end{lemma}

Our intention is now to settle the problem variationally. For that, we need to extend the function $f$ for $u<0$ and $u>L$. Given $\delta>0$, we extend $f$ to a function $\tilde{f} \in C^1(\mathbb{R})$ such that
\begin{equation*}
  \left \{
\begin{array}{lr}
\tilde{f}(u)=0 & \mbox{ if } u \geq L+\delta , \\
\tilde{f}(u)\leq 0& \mbox{ if } L\leq  u \leq L+\delta , \\
\tilde{f}(u)=f(u) & \mbox{ if } 0\leq u \leq L, \\
|\tilde{f}(u)|\leq 2\left(|f(0)|+1\right) & \mbox{ if } -\delta \leq u \leq 0,\\
\tilde{f}(u)=0 & \mbox{ if } u \leq -\delta.
\end{array} \right.
\end{equation*}
Notice that, since $F$ attains the maximum at $L$, $f'(L)=F''(L) \leq 0$. Thus, by combining with $f(L)=0$, we can extend $f$ such that  $\tilde{f}\leq 0$ in $[L,L+\delta]$ and this justifies the previous definition.

Accordingly, we define
\begin{equation*}
\tilde{F}(u) = \int_0^u \tilde{f}(s) \,ds+\frac{\sqrt{2}}{2}-1.
\end{equation*}

Note that there exists a sequence $\mu_n<L$, $\mu_n
\to L$, such that
\begin{equation}\label{e1.7}
\tilde F(\mu_n)> \tilde F(u),~\forall ~u <\mu_n.
\end{equation}
Indeed, $\mu_n$ can be constructed as follows. Define
$$m_n= \max\left\{\tilde F(x):\ x \in \left[0,L- \frac 1 n\right]\right\}, \ \mbox{and } \mu_n
= \min\left\{x \in \left[0,L- \frac 1 n\right]: \ \tilde F(x)=m_n\right\}.$$
By the definition of $\mu_n$, $\tilde{F}(\mu_n)=m_n > \tilde F(x)$ for all $x \in
[0,\mu_n)$. In addition, $m_n$ is increasing and $m_n\rightarrow \tilde{F}(L)=F(L)=0$. Clearly, $\mu_n < L$. We claim that $\mu_n \to L$. If not, there is a subsequence (still denoted by $\mu_n$) such that $\mu_n\rightarrow \mu<L$. Then $m_n=\tilde{F}(\mu_n)\rightarrow \tilde{F}(\mu)<0$, but $m_n\rightarrow 0$, which gives a contradiction.

Observe that
\begin{equation*}
\left|\tilde{F}(u)-\frac{\sqrt{2}}{2}+1\right|\leq 2(|f(0)|+1)\delta,~\forall ~u\leq 0.
\end{equation*}
We now fix $\delta>0$ (small) so that (note that $\tilde{F}(L)=0$)
\begin{equation}\label{e6.2}
 \frac{\sqrt{2}}{4}-1\leq \tilde{F}(u)\leq -\frac{1}{4}, ~~ \forall u \leq 0,
\end{equation}
\begin{equation}\label{e.close-u}
\tilde{F} (u)\geq -\frac{1}{8},~~~~\forall~~u\in [L-2\delta, L]
\end{equation}
and
\begin{equation} \label{cond-F}
0=\tilde{F}(L) > \tilde{F}(u), ~~ \forall u <L\,.
\end{equation}

With this truncation, our aim is to find radial solutions of the problem
\begin{equation}\label{eq-trunc}
\left \{
\begin{array}{ll} \mathrm{div}\left(\frac{\nabla u}{\sqrt{1+|\nabla u|^2}}\right) + \tilde{f}(u)=0~~& \mbox{ in } B_R,\\
u=0~~& \mbox{ on } \partial B_R,
\end{array} \right.
\end{equation}
for any $R>0$.

We use the calculus of variation to construct such solutions. Naturally, we should consider the energy functional:
\begin{equation*}
I_{R}(u)= \int_{B_R} \sqrt{1+|\nabla u|^2}- \tilde{F}(u).
\end{equation*}
However, it is well-known that with this energy functional, we can not obtain a minimizer in $W^{1,2}$ or $W^{1,1}$. Instead, the minimizer belongs to $BV(B_R)$ in general. Thus, we consider a perturbation of the above energy functional motivated by \cite[P. 163, Proof of Theorem 3.4]{MR526154}. For $\varepsilon>0$, define the approximation energy functional:
\begin{equation*}
  \begin{aligned}
&I_{R,\varepsilon}:\{u\in H^1_0(B_R): u \mbox{ is radially symmetric}\} \to \R,\\
&I_{R,\varepsilon}(u)= \int_{B_R} \sqrt{1+|\nabla u|^2}+\frac{\varepsilon}{2}|\nabla u|^2  - \tilde{F}(u).\\
  \end{aligned}
\end{equation*}
By the classical theory of calculus of variation and regularity theory, we have the:
\begin{lemma}\label{variational}
For any $R,\varepsilon>0$, the functional $I_{R,\varepsilon}$ attains its minimum at some function $u_{R,\varepsilon}$. Moreover, $u_{R,\varepsilon}\in C^{2,\alpha}(\bar{B}_R)$ for any $0<\alpha<1$ and
\begin{equation}\label{e1.4}
\begin{cases}
\mathrm{div}\left(\frac{\nabla u_{R,\varepsilon}}{\sqrt{1+|\nabla u_{R,\varepsilon}|^2}}\right)+\varepsilon \Delta u_{R,\varepsilon}+\tilde f(u_{R,\varepsilon})=0\quad & \mbox{in}~~B_R,\\
u_{R,\varepsilon}=0\quad & \mbox{on}~~\partial B_R.
\end{cases}
\end{equation}
\end{lemma}

In the following, we prove the existence of a minimizer of $I_R$ based on $u_{R,\varepsilon}$. First, we show the:
\begin{lemma}\label{mp} Let $u_{R,\varepsilon}$ be a solution of
\cref{e1.4}. Then
\begin{equation} \label{u-trunc}
-\delta<  u_{R,\varepsilon}< L ~~\mbox{ in}~~B_R.
\end{equation}
\end{lemma}

\proof First let us show that $u_{R,\varepsilon}(z)\leq L$ for any $z \in B_R$. Otherwise, assume that $\max u_{R,\varepsilon}= u_{R,\varepsilon}(z_0)>L$. Let
$\Omega=\{z \in B_R:\ u_{R,\varepsilon}(z) >L\}$. Then by the definition of $\tilde{f}$,
\begin{equation*}
\begin{cases}
\mathrm{div}\left(\frac{\nabla u_{R,\varepsilon}}{\sqrt{1+|\nabla u_{R,\varepsilon}|^2}}\right)+\varepsilon \Delta u_{R,\varepsilon} =-\tilde{f}(u_{R,\varepsilon}) \geq 0\quad & \mbox{in}~~\Omega,\\
u_{R,\varepsilon}=L\quad & \mbox{on}~~\partial \Omega.
\end{cases}
\end{equation*}
However, $u_{R,\varepsilon}$ attains its maximum in its interior, which is
impossible. In the same way we can prove that $u_{R,\varepsilon} \geq - \delta$.

We now show the strict inequality. Otherwise, assume that $\max u_{R,\varepsilon}
= L$. Observe also that the constant function $L$ is also a solution. Therefore both solutions are in contact, and this is in contradiction with the strong comparison principle. \proofend

The key for passing $\varepsilon \to 0$ is the following uniform gradient estimate, which relies heavily on $\inf \tilde{F}>-1$.

\begin{lemma}\label{variational-2}
Let $u_{R,\varepsilon}$ be as in \Cref{variational}. Then
\begin{equation}\label{e1.3}
|\nabla u_{R,\varepsilon}|^2\leq \frac{1}{(\inf \tilde F+1)^2}-1 \quad\mbox{in}~~B_R.
\end{equation}
\end{lemma}
\proof Since $u_{R,\varepsilon}$ is radial symmetric, \cref{e1.4} can be written as
\begin{equation}\label{e1.5}
u_{R,\varepsilon}''\left(1+u_{R,\varepsilon}'^2\right)^{-3/2}
+\frac{1}{r}u_{R,\varepsilon}'\left(1+u_{R,\varepsilon}'^2\right)^{-1/2}
+\varepsilon u_{R,\varepsilon}''+\frac{\varepsilon}{r}u_{R,\varepsilon}'+\tilde{f}(u_{R,\varepsilon})=0, ~\forall ~0<r<R,
\end{equation}
where $u_{R,\varepsilon}'$ denotes the derivative with respect to $r=|x|$. Let
\begin{equation*}
\mathcal{H}_\epsilon(r):=\tilde{F}(u_{R,\varepsilon})-\left(1+u_{R,\varepsilon}'^2\right)^{-1/2}+\frac{\varepsilon}{2} u_{R,\varepsilon}'^2.
\end{equation*}
With the aid of \cref{e1.5}, we have
\begin{equation*}
\mathcal{H}_\epsilon'=-\frac{1}{r}u_{R,\varepsilon}'^2\left(1+u_{R,\varepsilon}'^2\right)^{-1/2}
-\frac{\varepsilon}{r}u_{R,\varepsilon}'^2\leq 0.
\end{equation*}
Note that
\begin{equation*}
\mathcal{H}_\epsilon(0)=\tilde{F}(u_{R,\varepsilon}(0))-1<\tilde{F}(L)-1=-1.
\end{equation*}
Hence,
\begin{equation*}
\tilde{F}(u_{R,\varepsilon})-\left(1+u_{R,\varepsilon}'^2\right)^{-1/2}+\frac{\varepsilon}{2} u_{R,\varepsilon}'^2\leq -1.
\end{equation*}
From \cref{e6.1} and \cref{e6.2}, $\inf \tilde{F}>-1$. Then
\begin{equation*}
|u_{R,\varepsilon}'|^2\leq \frac{1}{(\inf \tilde F+1)^2}-1<+\infty
\end{equation*}
that is inequality \eqref{e1.3}.
\proofend

Once we have the uniform gradient bound for $u_{R,\varepsilon}$, \cref{e1.4} is uniformly elliptic and the ellipticity constants are independent of $\varepsilon$. Therefore, we have the:
\begin{lemma}\label{variational-3}
Let $u_{R,\varepsilon}$ be as in \Cref{variational}. Then for any $0<\alpha<1$,
\begin{equation}\label{e6.4}
\|u_{R,\varepsilon}\|_{C^{2,\alpha}(\bar{B}_R)}\leq C,
\end{equation}
where $C$ depends only on $n, \alpha, L$ and $\|f\|_{C^{0,1}}$.
\end{lemma}

Now, we can show that the minimizer of $I_R$ exists and is smooth.
\begin{lemma}\label{variational-4}
For any $R>0$, the functional $I_{R}$ attains its minimum at some function $u_{R}\in H^{1}_0(B_R)$ being radially symmetric. Moreover, $u_{R}\in C^{2,\alpha}(\bar{B}_R)$ for any $0<\alpha<1$ and
\begin{equation*}
\begin{cases}
\mathrm{div}\left(\frac{\nabla u_{R}}{\sqrt{1+|\nabla u_{R}|^2}}\right)+\tilde f(u_{R})=0\quad & \mbox{in}~~B_R;\\
u_{R}=0\quad & \mbox{on}~~\partial B_R.
\end{cases}
\end{equation*}
In addition, we have
\begin{equation}\label{e1.6}
-\delta<u_{R}< L~~\mbox{ in}~~ B_R
\end{equation}
and
\begin{equation}\label{e6.3}
|\nabla u_{R}|^2\leq \frac{1}{(\inf \tilde F+1)^2}-1 \quad\mbox{in}~~B_R,\quad
\|u_{R}\|_{C^{2,\alpha}(\bar{B}_R)}\leq C,
\end{equation}
where $C$ depends only on $n, \alpha, L$ and $\|f\|_{C^{0,1}}$.
\end{lemma}
\proof
By \Cref{variational-3}, up to a subsequence, there exists $u_R\in C^{2,\alpha}(\bar{B}_R)$ for any $0<\alpha<1$ such that $u_{R,\varepsilon}\to u_R$ in $C^{2}(\bar{B}_R)$. Then \cref{e6.3} holds by \cref{e1.3} and \cref{e6.4}. In addition, since for any $v\in H^1_0(B_R)$ radially symmetric
\begin{equation*}
I_{R,\varepsilon}(u_{R,\varepsilon})\leq I_{R,\varepsilon}(v),
\end{equation*}
then when $\varepsilon\to 0$ we obtain
\begin{equation*}
I_{R}(u_{R})\leq I_{R}(v).
\end{equation*}
Therefore, $u_R$ is a minimizer of $I_R$. Finally, the proof of \cref{e1.6} is similar to that of \cref{u-trunc} and we omit it.\proofend

In the rest of this section, we show that $u_R$ is positive for $R$ large enough. First, let us give some energy estimates on $u_R$. In what follows, we denote by $A_{R_1,R_2}(p)$ the annulus with center $p$ and radii
$R_1<R_2$. We have the:

\medskip

\begin{lemma} \label{lem-energy} There exists $C>0$
independent of $R$ so that:
\begin{equation} \label{estim-energy}
\pi R^2 \leq I_R(u_R)\leq \pi R^2 + C R, \end{equation}
and
\begin{equation} \label{estim-F}
- C R \leq \int_{B_R}\tF(u_R) \leq 0.
\end{equation}
\end{lemma}

\proof Taking into account \cref{cond-F} and \cref{e1.6}, we have immediately the second inequality of \cref{estim-F} and also
$$I_R(u_R) =
\int_{B_R}  \sqrt{1+|\nabla u_R|^2} - \tilde{F}(u_R) \geq
\int_{B_R} (1-\tilde{F}(u_R)) \geq \pi R^2,$$
which is the first inequality of \cref{estim-energy}.

To prove the second inequality of \cref{estim-energy}, let us define $\phi_R \in
H_0^1(B_R)$ by:
\begin{equation*}
  \phi_R(z) =
  \begin{cases}
   L & |z| \leq R-1, \\
  L(R-|z|) & |z| \in [R-1, R].
  \end{cases}
\end{equation*}
We now estimate $I_R(\phi_R)$. The gradient term can be estimated
as:
$$ \int_{B_R}  \sqrt{1+|\nabla \phi_R|^2}\leq \pi R^2+ \int_{B_R} |\nabla \phi_R|= \pi R^2+2 \pi \int_{R-1}^R
\phi_R'(r) r \, dr \leq \pi R^2+C R.$$
In order to estimate the term $\int_{B_R} \tilde{F}(\phi_R)$,
we split it into two terms:
\begin{equation*}
\int_{B_{R-1}} \tilde{F}(\phi_R)=0,\quad \int_{A_{R-1,R}(0)} \tilde{F}(\phi_R) \geq - C R
\end{equation*}
for a certain positive constant $C$.
In the last estimate we have just used the boundedness of
$\tilde{F}$.
The above estimates imply
\begin{equation*}
I_R(\phi_R) \leq \pi R^2 + C R.
\end{equation*}
Since $I_R(u_R) \leq I_R(\phi_R)$, we have the second inequality of
\cref{estim-energy}.
Finally,
$$ - \int_{B_R} \tilde{F}(u_R) = I_R(u_R)-\int_{B_R}  \sqrt{1+|\nabla u_R|^2} \leq C R,$$
and the first inequality of \cref{estim-F} follows.
\proofend

Next lemma is devoted to show the asymptotic behavior of
$u_R$.

\medskip

\begin{lemma} \label{prop} The following assertions hold:

\begin{enumerate}

\item[a)] For any fixed $0<\rho <L$, there exists $C=C_\rho$
independent of $R$ such that:
$$\Omega_\rho = \{ z \in B_R: u_R(z) < \rho\} \subset A_{R-C_\rho,
R}(0).$$

\item[b)] There exists $R_0>0$ such that $u_R$ is positive for $R
\geq R_0$.

\end{enumerate}
\end{lemma}

\proof The proof of a) will be made in two steps. In the first one we show that the volume of $\Omega_\rho$ is controlled by a linear growth in $R$. In the second we will show that $\Omega_\rho$ is connected. Such two facts will allow us to conclude that the $\Omega_\rho$ stays in a certain annulus.

\medskip

{\bf Step 1.} We show that for any fixed $\rho <L$, there exists $C=C_\rho$
independent of $R$ so that $|\Omega_\rho| \leq C_\rho R$. Indeed,
\begin{equation*}
 \int_{B_R \setminus \Omega_\rho} \tF(u_R) \leq 0,\quad
 \int_{\Omega_\rho} \tF(u_R) \leq \max\{\tF(x): x < \rho\}
|\Omega_\rho| = -\varepsilon |\Omega_\rho|,
\end{equation*}
where
$\varepsilon:= -\max\{\tF(x): x < \rho\}>0$ by
\cref{cond-F}.
Adding both terms, we get:
$$ \int_{B_R} \tilde{F}(u_R) \leq- \varepsilon
|\Omega_\rho|.$$ Then, from \cref{estim-F} we obtain
\[
|\Omega_\rho| \leq \frac{C}{\varepsilon} R\,
\]
and the conclusion follows with $C_\rho = C/\varepsilon$.

\medskip

{\bf Step 2.} Let us show that there exists $\mu$ such that $\Omega_\rho \subset \Omega_\mu$, and $\Omega_\mu$ is connected.

Take
\[
\mu:=\min \{\mu_n\, , \, \mu_n \geq \rho\}\,,
\]
where $\{\mu_n\}$ is the sequence in \cref{e1.7}. We claim now that $\Omega_{\mu}$ is
connected. Observe that $\Omega_{\mu}$ always has a connected
component touching the boundary $\partial B_R$, since $u_R = 0$ on $\partial B_R$. Suppose by
contradiction that it has an interior connected component too,
denoted by $U$. Then, $u_R(z) < \mu$ for $z \in U$ and $u_R(z)=
\mu$ if $z \in
\partial U$.

Define
$$ v(z)= \left \{ \begin{array}{lr} u_R(z) & z \notin U, \\ \mu & z
\in U. \end{array} \right. $$

Clearly, $v \in H_0^1(B_R)$ is radial symmetric and $\int_{U} \sqrt{1+|\nabla u_R|^2} >
\int_{U} \sqrt{1+|\nabla v|^2}$. Moreover, taking into account
\cref{e1.7},
$$\int_U \tF(u_R) < \int_U \tF(\mu) = \int_U \tF(v).$$
Therefore $I_R(v) < I_R(u_R)$, a contradiction since $u_R$ is a minimizer. So $\Omega_\mu$ is connected, and clearly $\Omega_\rho \subset \Omega_\mu$.

\medskip

From Steps 1 and 2 we have that $\Omega_\rho \subset \Omega_\mu \subset A_{R-C_\rho,
R}(0)$, concluding the proof.

\medskip

We now turn our attention to assertion b). Suppose that there exists $r_0 \in [0,R)$ with
$u_R(r_0)=-\delta_R \leq 0$, $u_R'(r_0)=0$. By \Cref{mp},
$\delta_R \in [0, \delta)$ and by a) we have that $r_0 \in
(R-C, R)$ for some positive $C>0$ independent of $R$.

Define $v(z)= u_R(z) + \delta_R$, which is a solution of the
problem:
\begin{equation*}
  \begin{cases}
  \mathrm{div}\left(\frac{\nabla v}{\sqrt{1+|\nabla v|^2}}\right) + g(v)=0 & \mbox{ in }B_{r_0}, \\
  v=0 & \mbox{ on } \partial B_{r_0}, \\
  \frac{\partial v}{\partial \nu}=0 & \mbox{ on } \partial B_{r_0}
  \end{cases}
\end{equation*}
where $g(t)= \tf(t- \delta_R)$. We now apply the Pohozaev's identity
(see \Cref{le1.1}) to the previous problem, to obtain that
\begin{equation} \label{cip}
\int_{B_{r_0}} G(v)-\frac{1}{2}\int_{B_{r_0}} \psi(|\nabla v|^2)=0,
\end{equation}
with
\begin{equation*}
G(t)= \int_0^t g(s) \, ds = \tF(t-\delta_R) -\tF(-\delta_R)
\end{equation*}
and
\begin{equation*}
\varphi(s)=2\left(\sqrt{1+s}-1\right),~
\psi(s)=2\left(\sqrt{1+s}-1-\frac{1}{2}\frac{s}{\sqrt{1+s}}\right).
\end{equation*}

We will show now that this is impossible if $R$ is sufficiently large.
Indeed, take $\Omega_\mu$ the set defined in Step 2. Then,
$$ \int_{B_{r_0}\setminus \Omega_{\mu}} G(v) = \int_{B_{r_0}\setminus
\Omega_{\mu}} \tF(u_R) - \tF(-\delta_R).$$
Noting that, and since $\mu$ is close to $L$, by \cref{e6.2} and \cref{e.close-u} we have
\begin{equation*}
\tF(u_R) - \tF(-\delta_R)\geq\frac{1}{8}\quad\mbox{ in }~~B_{r_0}\setminus \Omega_{\mu}\,.
\end{equation*}
Then
$$ \int_{B_{r_0}\setminus \Omega_{\mu}} G(v) \geq \frac{1}{8} | B_{r_0} \setminus \Omega_{\mu}| \geq c R^2.$$
In addition,
$$ \left |\int_{\Omega_\mu} G(v) \right | \leq \int_{A(0;R-C,R)}
|G(v)| \leq  CR.$$ Moreover, by \cref{estim-energy} and \cref{estim-F},
\begin{equation*}
  \begin{aligned}
\frac{1}{2}\int_{B_{r_0}} \psi(|\nabla v|^2)
=&\int_{B_{r_0}} \left(\sqrt{1+|\nabla v|^2}-1-\frac{1}{2}\frac{|\nabla v|^2}{\sqrt{1+|\nabla v|^2}}\right)\\
\leq &\int_{B_{R}(O)} \left(\sqrt{1+|\nabla u_R|^2}-1\right)\\
=&I(u_R)+\int_{B_{R}(O)} \tilde{F}(u_R)-\pi R^2\\
\leq& CR.
  \end{aligned}
\end{equation*}
Hence \cref{cip} cannot hold for large $R$. \proofend

\noindent \textbf{Proof of \Cref{th6.1}.} By combining \Cref{variational-4} and b) in \Cref{prop}, we obtain \Cref{th6.1}.~\qed~\\

\section{Proof of the main result: \Cref{Tmain}}\label{secM}
Before we prove our main theorem, we first show that the solutions in balls constructed in \Cref{th6.1} converge, in a suitable way, to a parallel solution in a half-space. Such result is valid in $\R^n$.
\begin{proposition} \label{teo-ball}
Let $u_R$ be the positive radial solution in \Cref{th6.1} for $R>R_0$. Let $p_R = (0,...,0,R) \in \R^n$. Then $v_R(x):= u_R(x-p_R)$ converges locally in compact sets of $\bar{H}:=\{x:=(x_1, ..., x_{n-1}, x_n) \, |\, x_n > 0\}$ to
\[
\Phi(x_1, ..., x_{n-1}, x_n) := \varphi(x_n)\,,
\]
where $\varphi$ is the positive bounded solution of \cref{eq.one}.
\end{proposition}
In order to simplify the notation, we will prove the result in dimension $n=2$, and we will denote as $x,y$ the two variables of $\R^2$. So $p_R=(0,R)$. Take a sequence $R_m \to +\infty$, $m \in \mathbb{N}$, $v_m(x,y)=u_{R_m}((x,y)-(0,R_m))$, $(x,y)\in  B_{R_m}(0,R_m)$. For any compact set $\Omega$ of $\bar{H}$, we have $\Omega\subset B_{R_m}(0,R_m)$ for $R_m$ large enough. By \Cref{th6.1},
\begin{equation*}
\|v_{m}\|_{C^{2,\alpha}(\bar\Omega)}\leq \|v_{m}\|_{C^{2,\alpha}(\bar{B}_{R_m}(0,R_m))}\leq C,
\end{equation*}
where $C$ is independent of $m$. As a consequence, $v_m$ converges (up to a subsequence) to a
solution $v$ of
\begin{equation*}
\begin{aligned}
\mathrm{div}\left(\frac{\nabla v}{\sqrt{1+|\nabla v|^2}}\right)+f(v)&=0\quad\mbox{in}~~\bar{H}.
\end{aligned}
\end{equation*}

Next, we show that $v$ is parallel. Take $p=(x, y) \in \bar{H}$ and $R_m$ large enough. We denote by $\rho_{m}$ its distance to $(0,R_m)$, that is, $\rho_{m}= \sqrt{x^2 +
(R_m-y)^2}$. Since $u_{R_m}$ is radially symmetric, then
$v_m(p)=v_m(0,R_m-\rho_{m})$. Observe now that $R_m-\rho_{m}
\to y$. Therefore $v_m(p) \to v(0,y)$, which is independent of
$x$. That is, $v$ is parallel and we will write $v(t)$ instead of $v(0,t)$ in the following argument.

We now prove that $v(0)=0$. Given $t>0$, by $v_m(0,0)=0$ and $\|v_{m}\|_{C^{2,\alpha}(\bar{B}_{R_m}(0,R_m))}\leq C$, we have
\begin{equation*}
|v_m(0,t)|=|v_m(0,t)-v_m(0,0)|\leq t\cdot \|Dv_{m}\|_{L^{\infty}(\bar{B}_{R_m}(0,R_m))}\leq Ct.
\end{equation*}
By taking $m\to \infty$,
\begin{equation*}
|v(t)|\leq Ct.
\end{equation*}
That is, $v(0)=0$.

Next, we prove
\begin{equation}\label{e7.1}
\lim_{t\to + \infty} v(t)=L.
\end{equation}
By \Cref{th6.1},
\begin{equation*}
v(t)=\lim_{m\to \infty} v_m(0,t) \leq L,~\forall ~t>0.
\end{equation*}
In addition, for any $0<\rho<L$, by a) in \Cref{prop}, there exists $C>0$ (depending on $\rho$) such that
\begin{equation*}
v_m(0,t) \geq \rho,~\forall ~t \in (C, 2 R_m-C).
\end{equation*}
As a consequence,
\begin{equation*}
v(t) \geq \rho,~\forall ~t \in (C, +\infty),
\end{equation*}
which implies \cref{e7.1}.

Finally, we will show that $v=\varphi$ given in \cref{eq.one}. Indeed, since $v$ is a parallel solution,
\begin{equation*}
\mathcal{H}(t):=F(v(t))-\frac{1}{\sqrt{1+|v'(t)|^2}}\equiv C
\end{equation*}
for some constant $C$ (see the proof of \Cref{th1.2}). Then with the aid of \cref{e6.1} and \cref{e7.1},
\begin{equation*}
C=\lim_{t\to \infty} \left(F(v(t))-\frac{1}{\sqrt{1+|v'(t)|^2}}\right)=-\lim_{t\to \infty}\frac{1}{\sqrt{1+|v'(t)|^2}}.
\end{equation*}
Thus, the limit of $v'(t)$ exists and must be $0$ (otherwise it would contradict the fact that $v$ is bounded). Hence, $C=-1$. Therefore, by \cref{e6.1} again,
\begin{equation*}
-1=F(0)-\frac{1}{\sqrt{1+|v'(0)|^2}}=\frac{\sqrt{2}}{2}-1-\frac{1}{\sqrt{1+|v'(0)|^2}}.
\end{equation*}
That is, $v'(0)=1$. In conclusion, $v=\varphi$ given in \cref{eq.one}.

Therefore we have proved the convergence of an adequate
subsequence. The uniqueness of the limit implies that actually the
whole $v_R$ converges for $R \to +\infty$. \proofend

\medskip

\noindent\textbf{Proof of \Cref{Tmain}.} According to \Cref{sec2}, we know that either $\Omega$ is a half-plane and $u$ is parallel, or the curvature of $\partial \Omega$ is strictly negative. Suppose that the latter holds. That is, $\Omega$ contains a
half-plane $\bar{H}$ internally tangent to $\partial \Omega$. We can suppose that $\bar{H}$ is the half-plane $\{y>0\}$ and the interior tangent point with $\partial \Omega$ is the origin.

Let $R$ large enough and consider the solution $u_R$ given by \Cref{th6.1}. Since the parallel solution $u_\infty$ is obtained as limit of a sequence $u_n$ of translations of the function $u$ in $\Omega$, we get that there exists a point $p \in \Omega$ such that the ball $B_R(p)$ is contained in $\Omega$ and the graph of the function $u_R$ defined in $B_R(p)$ stays under the graph of the function $u$ (note that $u_R$ has been defined in $B_R$, but by a translation in $\R^2$ can be defined in $B_R(p)$).

Next, we move the ball $B_R(p)$ inside $\Omega$ till it reaches the position of the ball $B_R(q)$ with $q=(0,R)$.
Observe that the graph of the function $u_R$, during the motion, cannot touch the graph of the function $u$ by the maximum
principle.

Consider now the function $v_R$ given by the \Cref{teo-ball} and its limit for $R \to +\infty$, that is clearly equal to $u_\infty(x,y):=\varphi(y) \leq u(x,y)$ for all $(x,y) \in \bar{H}$, where $\varphi$ is the positive bounded solution of \eqref{eq.one}. Moreover, the normal derivative of the functions $u$ and $u_\infty$ is the same at the origin, and by the maximum principle we get
\[
u=u_\infty\,.
\]
This shows that $\Omega=\bar{H}$, which contradicts the fact that the curvature of $\partial \Omega$ is strictly negative.

So $\Omega$ is a half-plane and $u$ is parallel.~\qed~\\

\section{Gradient estimate and proof of Theorems \ref{prop2.1} and \ref{cor_cap}}\label{sec_boundedDU}
This section is devoted to the proof of the gradient estimate of \Cref{prop2.1}. From this and \Cref{Tmain} we will deduce \Cref{cor_cap}. \Cref{prop2.1} is valid in $\R^n$, so in this section we will consider problem \cref{e.R} in dimension $n$.
The gradient bound has been proved in \cite{MR4519145} for $f\equiv 0$ and in \cite{MR4271788} for $f\equiv $ const. In fact, their proofs are applicable to \cref{e.R} with some modifications. The proof depends heavily on the following result:
\begin{proposition}\label{P.psi0}
Let $u$ be as in \Cref{prop2.1} and $(\Sigma,g)$ be the graph of $u$ with the induced metric from the Euclidean one in $\mathbb{R}^{n+1}$. That is,
\begin{equation*}
  \Sigma=\{(u(x),x):x\in \Omega\},~~g_{ij}=\delta^{ij}+u_iu_j.
\end{equation*}
Suppose that $\tilde{\Sigma}\subset \Sigma$ is connected and open with its closure $\overline{\tilde{\Sigma}}\subset \Sigma$.

Then for every $p_0\in \tilde\Sigma$, there exists $\psi_1\in C^{\infty}(\tilde\Sigma)$ satisfying
  \begin{equation}\label{e.psi}
  \begin{cases}
  \psi_1(p_0)=1, & \\
  \psi_1>1 &\mbox{ on $\tilde\Sigma \backslash \{p_0\}$, }\\
  \psi_1(p)\rightarrow +\infty &\mbox{ as $d(p,p_0)\rightarrow \infty$, }\\
  \Delta_g \psi_1\leq\psi_1 &\mbox{ on $\tilde\Sigma$,}
  \end{cases}
\end{equation}
where $d$ denotes the distance function on $(\Sigma,g)$.
\end{proposition}

\proof Since $u$ is bounded and $f$ is Lipschitz continuous, $f(u)$ is bounded. In addition, $u$ is constant on $\partial \Omega$. By \cite[Proposition 3.2, Case 2]{MR4271788}, we have
\begin{equation}\label{e.Bgb}
  \liminf_{r\rightarrow \infty} \frac{\log |B^g_r (p_0)\cap \Sigma|_g}{r^2} <\infty,
\end{equation}
where $B^g_r (p_0)$ is the geodesic ball on $(\Sigma,g)$ with radius $r$ and center $p_0$, and $|\cdot|_g$ denotes its Riemannian volume.

Since $(\Sigma, g)$ is a graph on (the connected domain) $\Omega\subset \mathbb{R}^n$ and  $\overline{\tilde{\Sigma}}\subset \Sigma$, then bounded subset of $\tilde \Sigma$ have compact closure in $\Sigma$. So, we can use \cite[Proposition 3.10]{MR4271788} obtaining that there exists a connected, complete Riemannian manifold $(N,h)$ and an isometric embedding $\phi: \tilde{\Sigma}\rightarrow N$ such that for $r$ large enough,
\begin{equation*}
  |B^h_r (\phi(p_0))|_h\leq 2\,|B^g_r (p_0)\cap \tilde{\Sigma}|_g+5\pi.
\end{equation*}
By combining with \cref{e.Bgb},
\begin{equation*}
\liminf_{r\rightarrow \infty} \frac{\log |B^h_r (\phi(p_0))|_h}{r^2} <\infty,
\end{equation*}
which implies, by \cite[Theorem 9.1]{MR1659871}, that $(N,h)$ is stochastically complete (see again \cite{MR1659871} for the definition of a stochastically complete Riemannian manifold).
Hence, according to \cite[Proposition 3.7]{MR4271788}), there exists $\psi_0\in C^{\infty}(N)$ satisfying
\begin{equation*}
  \begin{cases}
  \psi_0(\phi(p_0))=1, & \\
  \psi_0>1 &\mbox{ on $N \backslash \{\phi(p_0)\}$, }\\
  \psi_0(\phi(p))\rightarrow +\infty &\mbox{ as $d_h(\phi(p),\phi(p_0))\rightarrow \infty$ on $N$,}\\
  \Delta_h \psi_0\leq\psi_0 &\mbox{ on $N$.}
  \end{cases}
\end{equation*}
Then the function $\psi_1:=\psi_0(\phi)\in C^{\infty}(\tilde{\Sigma})$ satisfies \cref{e.psi}.\qed~\\

Before proving \Cref{prop2.1}, we prepare some preliminaries. Let $\{\partial_i\}_{i=1,...,n}$ be a local coordinate frame on $(\Sigma, g)$. Again we will use the index summation convention introduced in \Cref{sec2}. Define
\begin{equation*}
  W=\sqrt{1+|\nabla u|^2}=\sqrt{1+u_ku_k}.
\end{equation*}
Thus,
\begin{equation*}
  W_i=(1+|\nabla u|^2)^{-\frac{1}{2}}u_k u_{ki}=\frac{u_k u_{ki}}{W}.
\end{equation*}
In addition, we have
\begin{equation*}
  \det g=1+|\nabla u|^2=W^2.
\end{equation*}
The inverse of $g$ is given by
\begin{equation*}
  g^{ij}=\delta^{ij}-\frac{u_i u_j}{1+|\nabla u|^2}=\delta^{ij}-\frac{u_i u_j}{W^2}.
\end{equation*}
Besides, we calculate
\begin{equation}\label{e.2.11}
  \nabla_g u=g^{ij}u_j \partial_i=\frac{u_i \partial_i}{W^2}
\end{equation}
and
\begin{equation}\label{e.2.9}
\begin{aligned}
\Delta_g u=&\frac{1}{W}(W g^{ij}u_j)_i \\
  =&\frac{1}{W}\left(W \left(\delta^{ij}-\frac{u_i u_j}{W^2}\right)u_j\right)_i \\
  =&\frac{1}{W}\left(W \left(u_i-\frac{u_i |\nabla u|^2}{W^2}\right)\right)_i \\
  =&\frac{1}{W}\left(\frac{u_i}{W}\right)_i=-\frac{f(u)}{W}.
\end{aligned}
\end{equation}

\medskip

Now we can give the:~\\

\noindent\textbf{Proof of \Cref{prop2.1}.} Define the following linear operator:
\begin{equation}\label{e.2.1}
  L\, (\Phi):=W^2 \mathrm{div}_g (W^{-2}\nabla_g \Phi)=\Delta_g \Phi-2\left\langle \frac{\nabla_g W}{W}, \nabla_g \Phi \right\rangle_g,~\forall ~\Phi \in C^{2}(\Sigma).
\end{equation}
For any $\eta \in C^{2}(\Sigma)$, by a direct calculation, we have
\begin{equation}\label{e.2.2}
  \begin{aligned}
    L(\eta W)&=\Delta_g (\eta W)-2\left \langle \frac{\nabla_g W}{W}, \nabla_g (\eta W) \right \rangle_g\\
    &=\eta \Delta_g W+W \Delta_g \eta+2\left \langle \nabla_g \eta, \nabla_g W \right\rangle_g
    -2\left \langle \frac{\nabla_g W}{W}, W\nabla_g \eta +\eta\nabla_g W \right\rangle_g\\
    &=\eta \Delta_g W+W \Delta_g \eta-2\left \langle \frac{\nabla_g W}{W}, \eta \nabla_g W \right\rangle_g\\
    &=\eta LW+W\Delta_g \eta.
  \end{aligned}
\end{equation}
In addition,
\begin{equation}\label{e.2.3}
  \begin{aligned}
    L(W)=&\Delta_g W-2\left \langle \frac{\nabla_g W}{W}, \nabla_g W \right\rangle_g\\
    =&\Delta_g W-2\frac{|\nabla_g W|_g^2}{W}\\
    =&\frac{1}{W}(W g^{ij} W_j)_i-\frac{2}{W}g^{ij}W_iW_j \\
    =&\frac{1}{W}\left( \left(\delta^{ij}-\frac{u_i u_j}{W^2} \right) u_l u_{lj} \right)_i
   -\frac{2}{W}
   \left(\delta^{ij}-\frac{u_i u_j}{W^2} \right)\frac{u_ku_{ki}u_lu_{lj}}{W^2}\\
   =&\frac{1}{W}\left(u_{li}u_{li}+u_lu_{lii}
   -W^{-2}
   (u_{ii}u_ju_lu_{lj}+u_iu_{ji}u_lu_{lj}+u_iu_ju_{li}u_{lj}+u_iu_ju_lu_{lij})\right.\\
   &\left.+4W^{-4}u_iu_ju_lu_{lj}u_ku_{ki}
   -2W^{-2}u_ku_{ki}u_lu_{li} \right).
  \end{aligned}
\end{equation}
The equation
\begin{equation*}
  \mathrm{div} \left(\frac{\nabla u}{\sqrt{1+|\nabla u|^2}}\right)+f(u) =0
\end{equation*}
can be rewritten as
\begin{equation}\label{e.2.7}
  W^{-1}u_{ii}-W^{-3}u_iu_ju_{ij}+f(u)=0.
\end{equation}
By taking the first derivative with respect to $x_l$ in \Cref{e.2.7}, we have
\begin{equation}\label{e.2.4}
\begin{aligned}
  &W^{-1}u_{iil}
  -W^{-3}(u_ku_{kl}u_{ii}+u_{il}u_ju_{ij}+u_iu_{jl}u_{ij}+u_iu_ju_{ijl})\\
  &+3W^{-5}u_ku_{kl}u_i u_j u_{ij}+f'(u) u_l=0\,.
\end{aligned}
\end{equation}
Then, multiplying by $u_l$ on both sides and taking the sum over $l$, we obtain:
\begin{equation}\label{e.2.5}
\begin{aligned}
  &W^{-1}u_{iil}u_l
  -W^{-3}(u_ku_lu_{kl}u_{ii}+u_lu_{il}u_ju_{ij}+u_iu_lu_{jl}u_{ij}+u_iu_ju_lu_{ijl})\\
  &+3W^{-5}u_ku_l u_{kl}u_i u_j u_{ij}+f'(u) |\nabla u|^2=0\,.
\end{aligned}
\end{equation}
Inserting \Cref{e.2.5} into \Cref{e.2.3} we obtain
\begin{equation}\label{e.2.6}
  \begin{aligned}
    L(W)&=\frac{1}{W}\left(u_{li}u_{li}-f'(u) |\nabla u|^2 W
   +W^{-4}u_iu_ju_lu_{lj}u_ku_{ki}
   -2W^{-2}u_ku_{ki}u_lu_{li} \right).
  \end{aligned}
\end{equation}

Now we will show that $L(W)\geq 0$. Let
\begin{equation}\label{e.2.14}
  I:= u_{li}u_{li}+W^{-4}u_iu_ju_lu_{lj}u_ku_{ki}-2W^{-2}u_ku_{ki}u_lu_{li}.
\end{equation}
Since $I$ is invariant under orthogonal transformations of the coordinate system, we can fix a point $x$ and assume that $\nabla u(x)=(u_1,0,0,...,0)$, that is $u_i=0$ ($i\geq 2$) at the point $x$. Then, looking at the point $x$, \Cref{e.2.14} can be written as
\begin{equation*}
  I=u_{li}u_{li}+\frac{u_1^4u_{11}^2}{(1+u_1^2)^2}-\frac{2u_1^2 u_{1i} u_{1i}}{1+u_1^2}.
\end{equation*}
Note that
\begin{equation*}
  u_{li}u_{li}=\sum_{l=1}^{n} u_{l1}^2+\sum_{l=1}^{n }\sum_{i=2}^{ n}u_{li}^2.
\end{equation*}
Thus,
\begin{equation*}
  \begin{aligned}
    u_{li}u_{li}-\frac{2u_1^2 u_{1i} u_{1i}}{1+u_1^2}=&\left(\sum_{l=1}^{ n}u_{l1}^2\right)\frac{1-u_1^2}{1+u_1^2}
    +\sum_{l=1}^{ n}\sum_{i=2}^{n}u_{li}^2\\
    =&\left(\sum_{l=2}^{ n}u_{l1}^2\right)\frac{1-u_1^2}{1+u_1^2}
    +u_{11}^2\frac{1-u_1^2}{1+u_1^2}+\sum_{i=2}^{ n} u_{1i}^2
    +\sum_{l=2}^{n}\sum_{i=2}^{ n}u_{li}^2
  \end{aligned}
\end{equation*}
and then
\begin{equation*}
  I=\left(\sum_{l=2}^{ n}u_{l1}^2\right)\frac{2}{1+u_1^2}
    +\frac{u_{11}^2}{(1+u_1^2)^2}
    +\sum_{l=2}^{n}\sum_{l=2}^{n}u_{li}^2.
\end{equation*}
Hence, we obtain $L(W)\geq 0$ since $f'(u)\leq 0$.

For $\delta>0$ to be specified later, let
\begin{equation}\label{e.2.8}
\psi=\psi_1^{-1},~~~~\eta=\psi e^{-2u}-\delta,~~~~z=\eta W,
\end{equation}
where $\psi_1$ satisfies \Cref{e.psi} for some $p_0\in \Sigma$ to be specified later. Then
\begin{equation*}
  \frac{\nabla_g \psi}{\psi}=-\frac{\nabla_g \psi_1}{\psi_1}
\end{equation*}
and
\begin{equation}\label{e.2.10}
  \begin{aligned}
    \Delta_g \psi=&2\psi_1^{-3}|\nabla_g \psi_1|_g^2-\psi_1^{-2}\Delta_g \psi_1\\
    =&\psi\left(2\frac{|\nabla_g \psi_1|_g^2}{\psi_1^{2}}-\frac{\Delta_g \psi_1}{\psi_1}\right)\\
    \geq & \psi \left(2\frac{|\nabla_g \psi|_g^2}{\psi^{2}}-1\right).
  \end{aligned}
\end{equation}
Then, by noting \Cref{e.2.11}, \Cref{e.2.9} and \Cref{e.2.10} and the Young's inequality, we compute
\begin{equation*}
  \begin{aligned}
    \Delta_g (\psi e^{-2u})
    =&e^{-2u} \Delta_g \psi+\psi \Delta_g (e^{-2u})
    +2\left \langle \nabla_g \psi, \nabla_g (e^{-2u})\right \rangle_g\\
    =&e^{-2u} \Delta_g \psi+\psi e^{-2u}(4|\nabla_g u|_g^2-2\Delta_g u)
    -4e^{-2u}\left \langle \nabla_g \psi, \nabla_g u \right \rangle_g\\
    \geq &\psi e^{-2u}\left(2\frac{|\nabla_g \psi|_g^2}{\psi^{2}}-1
    +4|\nabla_g u|_g^2+\frac{2f}{W}
    -4\left \langle \frac{\nabla_g \psi}{\psi}, \nabla_g u \right \rangle_g\right)\\
    \geq &\psi e^{-2u}\left(2|\nabla_g u|_g^2-1+\frac{2f}{W}\right)\\
    =&\psi e^{-2u}\left(1-\frac{2}{W^2}+\frac{2f}{W}\right).
  \end{aligned}
\end{equation*}
Hence,
\begin{equation*}
  \begin{aligned}
    L(z)=L(\eta W)
    =&\eta L(W)+W\Delta_g \eta
    =\eta L(W)+W\Delta_g (\psi e^{-2u})\\
    \geq & \eta L(W)+\psi e^{-2u}\left(1-\frac{2}{W^2}+\frac{2f}{W}\right).
  \end{aligned}
\end{equation*}
Thus, if
\begin{equation}\label{e.z}
\eta\geq 0,~~W\geq A_0:= 8\left(\|f(u)\|_{L^{\infty}(\Omega)}+1\right)
\end{equation}
we have $L(z)\geq 0$.

Let
\begin{equation*}
  A_1:=\sup_{p\in \Sigma} W(p) e^{-2u(p)}.
\end{equation*}
We claim
\begin{equation}\label{e.2.12}
  A_1\leq \sqrt{1+\kappa^2}+A_0+1~~\mbox{ in}~~\Sigma.
\end{equation}
We prove it by contradiction. Suppose that $A_1>\sqrt{1+\kappa^2}+A_0+1$ ($A_1$ may be $+\infty$). Take $p_0\in \Sigma$ such that
\begin{equation*}
  W(p_0) e^{-2u(p_0)}>\sqrt{1+\kappa^2}+A_0+\frac{1}{2}.
\end{equation*}
Then by taking $\delta$ small enough such that $W(p_0)\delta<1/2$, we have $z(p_0)>\sqrt{1+\kappa^2}+A_0$.
We consider
\begin{equation*}
  \Sigma':=\{p\in \Sigma:z(p)>\sqrt{1+\kappa^2}+A_0\}.
\end{equation*}
Let $\tilde{\Sigma}\subset \Sigma'$ be the connected component with $p_0\in \tilde{\Sigma}$. It is obvious that $\partial \Sigma\cap \partial \tilde{\Sigma}=\emptyset$.
By the properties of $\psi_1$ \Cref{e.psi}, $\psi(p)\rightarrow 0$ as $d(p,p_0)\to \infty$. Thus, $\tilde{\Sigma}$ must be bounded. In addition, since $z(p_0)>\sqrt{1+\kappa^2}+A_0$, $z$ must attain the maximum at an interior point of $\tilde{\Sigma}$. Since
\begin{equation*}
  z=\eta W>\sqrt{1+\kappa^2}+A_0>0~~\mbox{ in }~~\tilde{\Sigma},
\end{equation*}
we get $\eta>0$ in $\tilde{\Sigma}$. Additionally, by noting that $\psi e^{-2u}<1$, we have $\eta<1$. Then,
\begin{equation*}
  W>\left(\sqrt{1+\kappa^2}+A_0\right)/\eta>A_0~~\mbox{ in }~~\tilde{\Sigma}.
\end{equation*}
Hence \cref{e.z} holds and then $Lz\geq0$ in $\tilde{\Sigma}$. The maximum principle gives a contradiction.

Thus, \cref{e.2.12} holds. That is,
\begin{equation*}
\sqrt{1+|\nabla u|^2}=W\leq \left(\sqrt{1+\kappa^2}+A_0+1\right)e^{2\|u\|_{L^{\infty}(\Omega)}}~~\mbox{ in }~~\Sigma.
\end{equation*}
Therefore, \cref{e2.13} holds and this concludes the proof.~\qed

\medskip

Next, let's give the:

\medskip

\textbf{Proof of \Cref{cor_cap}.} Let $u$ be the solution of \cref{capil-0} with $\Omega$ diffeomorphic to the half-plane $\mathbb{R}_{+}^2$. By the comparison principle for unbounded domains (see \cite{MR991023}), $0<u< c_h$. Set $v=c_h-u$ and then $v$ is a bounded solution of \cref{capil}. Additionally, since $f(v)=-bv+bc_h$, we obtain that $f\in C^{\infty}(\mathbb{R})$ and $f'(v)=-b<0$. By the gradient estimate given by \Cref{prop2.1}, $|\nabla v|$ is bounded. Then the equation is uniformly elliptic. By the standard Schauder regularity \cite[Theorem 6.17]{MR1814364} and $f\in C^{\infty}(\mathbb{R})$, $v\in C^{\infty}(\Omega)$.

Moreover, from \cite[Corollary 2]{MR607989}, there exists a unique parallel solution of \cref{capil-0} (in \cite{MR607989} the capillary problem is considered with $b=1$, but for a general $b$ one can apply the normalization lemma \cite[Lemma 2]{MR607989}). 

Therefore, the assumptions of \Cref{co1.1} are all satisfied. Then $\Omega$ is a half-plane and $v$ is parallel. That is, $u$ is parallel and the infinite plate must be planar.~\qed~\\



%
%
\printbibliography

\end{document}